\documentclass{article}
\usepackage{amssymb,latexsym,amscd}

\newcommand{\se}[1]{{\section{#1}} {\setcounter{equation}{0}}}
\newtheorem{theorem}{Theorem}[section]
\newtheorem{lm}{Lemma}[section]
\newtheorem{prop}{Proposition}[section]

\def\k{{K\"{a}hler }}

\begin{document}
\hbadness=10000
\title{{\bf Generalized special Lagrangian fibration for Calabi-Yau hypersurfaces in toric varieties III:}\\
{\Large {\bf The smooth fibres}}}
\author{Wei-Dong Ruan\\
Department of Mathematics\\
University of Illinois at Chicago\\
Chicago, IL 60607\\}
\footnotetext{Partially supported by NSF Grant DMS-0104150.}
\maketitle
\begin{abstract}
In this paper we construct all smooth torus fibres of the generalized special Lagrangian torus fibrations for Calabi-Yau hypersurfaces in toric varieties near the large complex limit.
\end{abstract}
\se{Introduction}
This paper is a sequel of \cite{sl1,sl2}. The aim of this series of papers is to construct generalized special Lagrangian torus fibrations for Calabi-Yau hypersurfaces in toric varieties near the large complex limit.\\
 
Let $(P_\Delta,\omega)$ be a toric variety whose moment map image (with respect to the toric \k form $\omega$) is the real convex polyhedron $\Delta \subset M_{\mathbb{R}}$. Also assume that the anti-canonical class of $P_\Delta$ is represented by an integral reflexive convex polyhedron $\Delta_0\subset M$ and the unique interior point of $\Delta_0$ is the origin of $M$. Integral points $m \in \Delta_0$ correspond to holomorphic toric sections $s_m$ of the anti-canonical bundle. For the unique interior point $m_o$ of $\Delta_0$, $s_{m_o}$ is the section of the anti-canonical bundle that vanishes to first order along each toric divisor of $P_\Delta$.\\

Let $\{w_m\}_{m\in \Delta_0}$ be a strictly convex function on $\Delta_0$ such that $w_m> 0$ for $m\in \Delta_0\setminus \{m_o\}$ and $w_{m_o}\ll 0$. Define

\[
\tilde{s}_t = s_{m_o} + ts,\ \ s = \sum_{m\in \Delta_0\setminus \{m_o\}} a_m s_m,\  {\rm where}\ |a_m| = \tau^{w_m},\ {\rm for}\ m\in \Delta_0\setminus \{m_o\}.
\]

(As in the case of Fermat type quintic, in general, $\Delta_0$ need not contain all the integral $m$ in the real polyhedron spanned by $\Delta_0$.) Let $X_t = \{\tilde{s}_t^{-1}(0)\}$. Then $\{X_t\}$ is a 1-parameter family of Calabi-Yau hypersurfaces in $P_{\Delta}$. $X_0 = \{s_{m_o}^{-1}(0)\}$ is the so-called large complex limit. $X_t$ is said to be near the large complex limit if $\tau$ and $t$ are small and $t \leq \tau^{-w_{m_o}}$.\\

{\bf Remark:} In fact, for the discussions in this paper, the coefficients $\{a_m\}_{m\in \Delta_0\setminus \{m_o\}}$ do not need to have the above specific form, and ``near the large complex limit" can also be relaxed to only require that $t$ is small. Although in the above specific situation, the fibration constructed will have more interesting structure and is better understood.\\

Since $X_0$ is toric, the moment map induces the standard generalized special Lagrangian fibration $F_0: X_0 \rightarrow \partial \Delta$ with respect to the toric holomorphic volume form. In \cite{lag1, lag2, lag3, toric}, we constructed Lagrangian torus fibration for $X_t$ when $X_t$ is near the large complex limit, using the Hamiltonian-gradient flow to deform the fibration $F_0$ for $X_0$ symplectically to the desired Lagrangian fibration $F_t: X_t \rightarrow \partial \Delta$ for such $X_t$. The (topological) singular set of the fibration map $F_t$ is $C= X_t \cap {\rm Sing}(X_0)$, which is independent of $t$. The corresponding singular locus $\tilde{\Gamma} = F_0(C)$ is also independent of $t$. When $X_t$ is near the large complex limit, $\tilde{\Gamma} \subset \partial \Delta$ exhibits the amoeba structure that is a fattening of a graph $\Gamma \subset \partial \Delta$. It was conjectured in \cite{lag2} that the singular locus for the (generalized) special Lagrangian torus fibration should resemble the singular locus $\tilde{\Gamma}$ of the Lagrangian torus fibration.\\ 

\cite{sl1} initiated the program of constructing generalized special Lagrangian fibration for $X_t$ using similar ideas by deformation from the standard fibration $F_0$ of the large complex limit $X_0$.  Since $X_0$ is a singular variety and $F_0: X_0 \rightarrow \partial \Delta$ is a singular fibration, it is a singular deformation problem. Such singular deformation problems are usually rather delicate and difficult if at all possible to be solved. Even for the much softer symplectic singular deformation of Lagrangian fibrations we discussed in \cite{lag1, lag2, lag3, toric}, the technical difficulties are already quite delicate and daunting.\\

One crucial observation as pointed out in \cite{sl1} is that, due to the canonical nature of the generalized special Lagrangian fibration on a Calabi-Yau manifold when we fix the \k metric, one can construct the generalized special Lagrangian fibration over different parts of the base $\partial \Delta$ separately and they will automatically match on the overlaps. More precisely our strategy of construction is, for different parts of $X_t$, we construct explicit generalized special Lagrangian fibrations for smooth local models $Y_t$, which are small perturbations of the corresponding regions of $X_t$. Then we construct smooth deformation family connecting $Y_t$ with $X_t$. Through smooth deformation, we get the generalized special Lagrangian fibration for this part of $X_t$. Through this idea, we turn the singular deformation problem into smooth deformation problems, which are much easier to handle. This idea in \cite{sl1} is the fundamental idea that makes all the subsequent works possible. Another key ingredient of \cite{sl1} is the quantitative implicit function theorem. In \cite{sl1}, using a local model that dates back to \cite{HL}, we are able to construct monodromy representing generalized special Lagrangian torus fibration for Fermat type quintic Calabi-Yau hypersurfaces in $\mathbb{CP}^4$ near the large complex limit.\\ 

Recall from \cite{sl1} that a (generalized) special Lagrangian fibration of smooth torus over an open set $U\subset \partial \Delta$ is said to represent the monodromy if $\partial \Delta \setminus U$ is a fattening of $\tilde{\Gamma}$ that retracts to $\tilde{\Gamma}$. \\

To deal with more general situations than those in \cite{sl1}, it is necessary to deal with ``thin torus" fibres that have small circles in some directions and big circles in other directions. Although the deformation of ``thin torus" fibres is smooth deformation, the classical estimates usually depend on different scales of the ``thin torus" and are not very useful when the torus is thin. Such difficulty is in a way similar to the difficulty one might encounter in singular deformations. To get meaningful results, it is crucial to derive estimates that do not depend on the multi-scales of the ``thin torus" fibres. In \cite{sl2}, we developed such multi-scale estimates for ``thin torus" fibres. As application, we are able to construct monodromy representing generalized special Lagrangian torus fibration for the mirror of quintic Calabi-Yau hypersurfaces in $\mathbb{CP}^4$ near the large complex limit.\\ 

Further more, in \cite{sl2}, using ideas related to amoeba, we are able to use similar local model as in \cite{sl1} for certain region of $X_t$, where one of $a_ms_m$ term dominates. Combined with our multi-scale estimates, we are able to construct monodromy representing generalized special Lagrangian torus fibration for Calabi-Yau hypersurfaces in toric varieties, including general quintic Calabi-Yau hypersurfaces in $\mathbb{CP}^4$, near the large complex limit. Due to the limitation of the local models, this construction in \cite{sl2} is less satisfactory, because we have to require a more strict sense of near the large complex limit and the toric metrics have to be specifically chosen accordingly.\\

It will be ideal if one can construct maximal monodromy representing generalized special Lagrangian torus fibration, which amounts to a generalized special Lagrangian torus fibration over $\partial \Delta$ away from a small neighborhood of the singular locus $\tilde{\Gamma}$. The constructions of monodromy representing generalized special Lagrangian torus fibrations in \cite{sl1,sl2} are far from being maximal. Such maximal construction is essentially the construction of all the smooth fibres. In this paper, we will give construction of maximal monodromy representing generalized special Lagrangian torus fibration for general Calabi-Yau hypersurfaces near the large complex limit in toric varieties with arbitrary smooth toric metrics.\\

The main ingredient that makes this paper possible is the construction of explicit local models in section 2. It is still quite mysterious to me that such explicit local models actually exist for every part of $\partial \Delta \setminus \tilde{\Gamma}$, yet unlike in \cite{sl2}, these local models are so accurate that they enable us to prove the construction in the most general situations.\\

While the torus fibres in the local models of \cite{sl1, sl2} are flat, the thin torus fibres in the local models in this paper are not flat and the different scales are more mixed than in \cite{sl2}. Interestingly, the multi-scale estimates for such thin torus, although much more delicate and involved, will still work.\\ 

One major problem in constructing (generalized) special Lagrangian fibration is to ensure that a (generalized) special Lagrangian torus does not intersect with nearby torus, which is equivalent to that the deformation 1-forms do not have zeroes. In \cite{sl1,sl2}, the (generalized) special Lagrangian torus is a small perturbation of flat torus, and the deformation 1-forms are small perturbations of the flat 1-forms, therefore have no zeroes. In this paper, the (generalized) special Lagrangian torus is not a small perturbation of flat torus. Yet the deformation 1-forms are somehow still small perturbations of flat 1-forms, therefore have no zeroes. This is another indication of the magical properties of our explicit local models.\\   

In this paper, we essentially construct all the smooth fibres for the global generalized special Lagrangian fibration. The next step would be to fill in the singular fibres. Notice that the local models in section 2 also make sense for singular fibres in certain situations. In the sequels of this paper, we will discuss some singular fibres.\\

This paper is organized as follows. In section 2 we construct the explicit local models of generalized special Lagrangian torus fibrations. In section 3, we discuss the concept of (strongly) T-boundedness that is closely related to the toric properties and is crucial for our multi-scale estimates. In section 4, we work out the Darboux coordinates of the symplectic neighborhood of a special Lagrangian torus in the local model and their estimates. The corresponding discussions in \cite{sl1,sl2} are much easier and almost trivial in comparison. In section 5, we establish the estimates for the deformation from the local models discussed in section 2 to the actual Calabi-Yau hypersurfaces. In section 6, we apply the symplectic neighborhood estimate for the local model $(Y_t, \check{\omega}_t)$ in section 4 and the estimate for the coordinate under deformation in section 5 to work out the basic estimates necessary to apply the implicit function theorem to our deformation construction. In section 7, we clarify why individual generalized special Lagrangian torus we construct via deformation will together form a fibration in the corresponding region. In section 8, we prove the local uniqueness of generalized special Lagrangian torus under Hamiltonian deformation, which makes it possible for us to glue different pieces of fibrations together. In section 9, we apply the results from section 7, 8 and \cite{sl2} to construct the global maximal monodromy representing generalized special Lagrangian torus fibration for general Calabi-Yau hypersurfaces near the large complex limit in toric varieties.\\

{\bf Remark on notations:} $a\sim b$ means $|\log (a/b)|$ is bounded or equivalently $C^{-1}|a| \leq |b| \leq C|a|$ for a bounded constant $C>0$. We often use $(a_{ij})$ (or $(a_{ij})_{m\times n}$) to denote a ($m\times n$) matrix with the $ij$-th entry $a_{ij}$. We also use $(A)_{ij}$ to denote the $ij$-th entry of a matrix $A$. $\alpha$ has been used for two purposes, as a $1$-form or in the notation $C^\alpha$ for Horder space. $I$ is used as either identity matrix or an index set. It should be clear from the context which usage is intended.\\

\se{Examples of special Lagrangian}
In this section, we discuss some examples of smooth generalized special Lagrangian tori that will be the local models of our deformation argument.\\

Consider $\mathbb{C}^{n+1}$ with coordinate $\tilde{z} = (z_0, z'', z') = (\tilde{z}'',z') = (z_0,\cdots,z_n)$ and toric metric with \k form $\omega = \partial\bar{\partial}(|\tilde{z}''|_\lambda^2 + \rho(z'))$, where $z'' = (z_1,\cdots,z_l)$, $z'= (z_{l+1},\cdots,z_n)$, $\displaystyle |\tilde{z}''|_\lambda^2 = \sum_{k\in I''} \lambda_k(z') |z_k|^2$, $I'' = \{0,\cdots,l\}$ and $I' = \{l+1,\cdots,n\}$. \\

Let $Y_t = \{\tilde{z}\in \mathbb{C}^{n+1}|z_0\cdots z_n = tp(z')\}$ be a family of hypersurfaces in $\mathbb{C}^{n+1}$.  $z = (z'',z') = (z_1,\cdots,z_n)$ can be used as coordinate of $Y_t$.\\

For any $r=(r_{l+1},\cdots,r_n) \in \mathbb{R}_{\geq 0}^{n-l}$ and $c=(c_0,\cdots,c_l) \in \mathbb{R}_{\geq 0}^{l+1}$, one can consider 

\[
S_{r,c} = \{\tilde{z} = (\tilde{z}'',z')\in \mathbb{C}^{n+1}||z_j|=r_j, j\in I'; \lambda_k|z_k|^2 = c_k + \eta(z'), k\in I''\},
\]

where $\eta(z')$ is a real valued function satisfying

\[
\prod_{k\in I''}(c_k + \eta(z')) = \kappa(z') = \tilde{\Lambda}(z')|tp(z')|^2,\ \ \tilde{\Lambda}(z') = \Lambda(z')R^{-1}(z'),
\]
\[
\Lambda = \prod_{k\in I''}\lambda_k(z'),\ \ R(z') = \prod_{j\in I'}r_j^2,
\]
with fixed $|t|$. Without loss of generality, we assume $0=c_0 \leq \cdots \leq c_l$. In particular, when $c=0$,

\[
S_{r,0} = \{\tilde{z} = (\tilde{z}'',z')\in \mathbb{C}^{n+1}||z_{j}|=r_j, j\in I'; \lambda_k|z_k|^2 = \tilde{\Lambda}^{\frac{1}{l}}|tp(z')|^{\frac{2}{l}}, k\in I''\}.
\]

It is convenient to introduce $\check{\eta}$ satisfying

\[
\prod_{k\in I''}(c_k + \check{\eta}) = \tilde{\Lambda}(z')|t|^2.
\]

Notice that $\tilde{\Lambda}(z')$ is actually a function of $r$. Hence $\check{\eta}$ is a function of ${r,c}$, which will be constant on $S_{r,c}$. It is easy to see that $\eta \leq C\check{\eta}$.\\

{\bf Remark:} $\eta$, $\kappa$, $S_{r,c}$ also depend on $|t|$, which we will take to be a fixed small constant throughout our discussion.\\

A region (in $\mathbb{C}^{n+1}$) is called {\bf normal}, if $|p(z')|\geq C>0$ and $\displaystyle \epsilon(\tilde{z}) = \max_{j\in I'} \left(\frac{|z_0|}{|z_j|}\right)$ is small in this region. (The notion of normal region will be refined in section 6 to suit our application. Such refinement is delayed to allow us to derive our estimates in more general setting before section 6.) It is easy to see that $\eta \sim \check{\eta}$ in a normal region. In a normal region near $S_{r,c}$, $|z_k|\sim \eta + c_k$ for $k\in I''$ and $\epsilon \sim \check{\eta}/r_{l+1}$. Without further mentioning, the discussions in this paper will always be within a normal region.\\

It is straightforward to derive

\begin{equation}
\label{be}
\frac{d \eta}{d \log \kappa} = \zeta\eta,\ \ {\rm where}\ \zeta = \left(\sum_{k\in I''}\eta(c_k + \eta)^{-1}\right)^{-1},\ 1\leq \zeta^{-1} \leq |I''|.
\end{equation}

Let $\mu$ satisfy

\[
\frac{d\mu}{d\log \kappa} = \eta.
\]

Then

\[
d\mu = \eta \sum_{k\in I''}\frac{d\eta}{c_k + \eta} = \sum_{k\in I''}d(\eta - c_k\log (c_k + \eta)).
\]

We may take

\[
\mu = (l+1)\eta - \sum_{k\in I''}c_k\log (c_k + \eta) = \eta(l+1 + \log \kappa) - \sum_{k\in I''}(c_k + \eta)\log (c_k + \eta).
\]

It is easy to check that

\begin{equation}
\label{bg}
\frac{\partial \eta}{\partial c_k} = - \frac{\zeta\eta}{c_k + \eta},\ \ \frac{\partial \mu}{\partial c_k} = - \log(c_k + \eta).
\end{equation}

Straightforward computation yields

\begin{equation}
\label{bf}
\partial\bar{\partial}\mu = \eta(\partial\bar{\partial}\log \Lambda + \zeta\partial\log \kappa\bar{\partial}\log \kappa).
\end{equation}

Assume that $z_k = r_k e^{i\theta_k}$. Let $\rho_{j} = \frac{\partial \rho}{\partial \log |z_j|^2}$, $\lambda_{k,j} = \frac{\partial \lambda_k}{\partial \log |z_j|^2}$, $\tilde{\Lambda}_j = \frac{\partial \tilde{\Lambda}}{\partial \log |z_j|^2}$ for $j\in I'$. We have\\
\begin{prop}
\label{bd}
\[
\omega - \partial\bar{\partial}\mu = -i\sum_{k\in I''} d(\lambda_k|z_k|^2 - \eta)\wedge d\theta_k -i \sum_{j\in I'} d\left(\rho_{j} + \sum_{k\in I''}\frac{\lambda_{k,j}}{\lambda_k}(\lambda_k|z_k|^2 - \eta)\right)\wedge d\theta_j.
\]
\end{prop}
{\bf Proof:} It is straightforward to derive that

\[
\omega = -i\sum_{k\in I''} d(\lambda_k|z_k|^2)\wedge d\theta_k -i \sum_{j\in I'} d\left(\rho_{j} + \sum_{k\in I''}\lambda_{k,j}|z_k|^2\right)\wedge d\theta_j,
\]

\[
\partial\bar{\partial}\mu =  -id\left(\frac{d\mu}{d\log \kappa}\right) \wedge d{\rm Arg}(tp) - i\sum_{j\in I'}d\left(\frac{d\mu}{d\log \kappa}\frac{\partial\log \tilde{\Lambda}}{\partial \log |z_j|^2}\right) \wedge d\theta_j
\]
\[
= -id\eta \wedge d{\rm Arg}(tp) - i\sum_{j\in I'}d\left(\eta\frac{\tilde{\Lambda}_j}{\tilde{\Lambda}}\right) \wedge d\theta_j.
\]

Subtract $\partial\bar{\partial}\mu$ from $\omega$ and use the facts

\[
d{\rm Arg}(tp) = \sum_{k=0}^n d\theta_k,\ {\rm and}\ \frac{\tilde{\Lambda}_j}{\tilde{\Lambda}} + 1 =\sum_{k\in I''}\frac{\lambda_{k,j}}{\lambda_k},
\] 

we get the desired formula.
\hfill\rule{2.1mm}{2.1mm}\\
\begin{prop}
\label{ba}
\[
\omega|_{S_{r,c}} = (\partial\bar{\partial}\mu)|_{S_{r,c}}.
\]
\end{prop}
{\bf Proof:} On $S_{r,c}$, $\rho_{k}$, $\lambda_{k,j}$, $(\lambda_k|z_k|^2 - \eta)$ are constants. By proposition \ref{bd},

\[
\omega|_{S_{r,c}} = (\partial\bar{\partial}\mu)|_{S_{r,c}}.
\]
\hfill\rule{2.1mm}{2.1mm}
\begin{theorem}
\label{bb}
When $|t|$ is small, $\check{\omega} = (\omega - \partial\bar{\partial}\mu)$ is a \k form on a normal region in $\mathbb{C}^{n+1}$. $S_{r,c}$ inside such region is a generalized special Lagrangian submanifold in $\mathbb{C}^{n+1}$ with respect to the symplectic form $\check{\omega}$ and the holomorphic volume form 
\[
\Omega = \prod_{k=0}^{n}\frac{dz_k}{z_k}.
\]
\end{theorem}
{\bf Proof:} Recall from the formula (\ref{bf})

\[
\partial\bar{\partial}\mu = \eta(\partial\bar{\partial}\log \Lambda + \zeta\partial\log \kappa\bar{\partial}\log \kappa).
\]

Since $|t|$ is constant in $\kappa$, 

\[
\partial\log \kappa = \partial(\log \tilde{\Lambda} + \log |p(z')|^2) = \partial(\log \Lambda + \log |p(z')|^2) - \sum_{j\in I'}\frac{dz_j}{z_j}.
\]

Since $|p|$ is bounded from below by a positive constant, $\partial(\log \Lambda + \log |p(z')|^2)$ is bounded. Consequently, $\partial\bar{\partial}\mu = O(\epsilon^2)$, which is small. Therefore, $\check{\omega}$ is a \k form on a normal region in $\mathbb{C}^{n+1}$.\\

With proposition \ref{ba}, to show that $S_{r,c}$ is a generalized special Lagrangian, we only need to show that $\Omega$ has constant phase when restricted to $S_{r,c}$. On $S_{r,c}$,

\[
\frac{dz_k}{z_k} = id\theta_k,\ {\rm for}\ k\in I';\ \frac{dz_k}{z_k} = id\theta_k + \frac{d\eta}{2(c_k + \eta)},\ {\rm for}\ k\in I''.
\]

When restricted to $S_{r,c}$, $|t|$ is a constant and $\eta$ only depends on $\theta'$.

\[
\Omega|_{S_{r,c}} = i^{n+1}\prod_{k=0}^{n}d\theta_k.
\]
\hfill\rule{2.1mm}{2.1mm}
\begin{theorem}
\label{bc}
When $|t|$ is small, $S_{r,c}\cap Y_t$ in a normal region is a generalized special Lagrangian submanifold in $Y_t$ with respect to the symplectic form $\check{\omega}_t = (\omega - \partial\bar{\partial}\mu)|_{Y_t}$ and the holomorphic volume form 
\[
\Omega_t = \prod_{k=1}^{n}\frac{dz_k}{z_k}.
\]
\end{theorem}
{\bf Proof:} We only need to show that $\Omega_t$ has constant phase when restricted to $S_{r,c}\cap Y_t$. The proof is almost the same as in theorem \ref{bb}. We have

\[
\Omega_t|_{S_{r,c}\cap Y_t} = i^{n}\prod_{k=1}^{n}d\theta_k.
\]
\hfill\rule{2.1mm}{2.1mm}\\
{\bf Remark:} Since $\check{\omega}_t$ depend on $c$, the family $\{S_{r,c}\cap Y_t\}$ generally does not form a generalized special Lagrangian fibration for $Y_t$ with respect to a fixed \k form.\\\\

\se{T-boundedness}
T-boundedness discussed in this section is closely related to toric properties. It will be useful in estimates throughout the paper. (In this section, $c$ is temporarily used to denote a constant, because $C$ sometimes is too easily confused as a matrix.)\\

Let $Z = {\rm Diag}(z_1,\cdots,z_n)$. A Hermitian matrix $G = (g_{i\bar{j}})_{n\times n}$ is called T-bounded if $ZGZ^{-1}$ is bounded. Notice that as a consequence of T-boundedness, $G$ and $\bar{Z}^{-1}G\bar{Z}$ are bounded.\\

\begin{lm}
\label{ca}
If $G$ is T-bounded and $\det (G) \geq c >0$, then $G^{-1}$ is T-bounded.
\end{lm}
{\bf Proof:} $G$ being bounded and $\det (G) \geq c >0$ imply that $G^{-1}$ is bounded. Since $\det (ZGZ^{-1}) = \det (G)$, $(ZGZ^{-1})^{-1} = ZG^{-1}Z^{-1}$ is bounded. Consequently, $G^{-1}$ is T-bounded.
\hfill\rule{2.1mm}{2.1mm}\\

{\bf Remark:} When applying this lemma, it is often convenient to verify the stronger conditions $G \geq cI$ ($c>0$) or $G^{-1}$ being bounded in place of $\det (G) \geq c >0$.\\

For our application, $G$ will be the metric matrix of a toric \k metric. More generally, we can consider a symplectic form that is a small perturbation of a \k form. Such symplectic form can be expressed as

\[
S = \hat{G} + \hat{H} = \left(\begin{array}{cc}H&G\\-G^t&-\bar{H}\end{array}\right),
\]

where

\[
\hat{G} = \left(\begin{array}{cc}0&G\\-G^t&0\end{array}\right),\ \hat{H} = \left(\begin{array}{cc}H&0\\0&-\bar{H}\end{array}\right).
\]

\[
\bar{\hat{Z}}S\hat{Z}^{-1} = \left(\begin{array}{cc}ZH\bar{Z}^{-1}&ZGZ^{-1}\\-\bar{Z}G^t\bar{Z}^{-1}&-\bar{Z}\bar{H}Z^{-1}\end{array}\right),\ {\rm where}\ \hat{Z} = \left(\begin{array}{cc}\bar{Z}&0\\0&Z\end{array}\right).
\]

\[
\hat{G}^{-1} = \left(\begin{array}{cc}0&-(G^{-1})^t\\G^{-1}&0\end{array}\right),\ \hat{Z}\hat{G}^{-1}\bar{\hat{Z}}^{-1} = \left(\begin{array}{cc}0&-\bar{Z}(G^{-1})^t\bar{Z}^{-1}\\ZG^{-1}Z^{-1}&0\end{array}\right). 
\]

The boundedness of $G$ and $ZGZ^{-1}$ is equivalent to the boundedness of $\hat{G}$ and $\bar{\hat{Z}}\hat{G}\hat{Z}^{-1}$. The boundedness of $G^{-1}$ and $ZG^{-1}Z^{-1}$ is equivalent to the boundedness of $\hat{G}^{-1}$ and $\hat{Z}\hat{G}^{-1}\bar{\hat{Z}}^{-1}$. $S$ is called T-bounded if $\bar{\hat{Z}}S\hat{Z}^{-1}$ is bounded. $S^{-1}$ is called T-bounded if $\hat{Z}S^{-1}\bar{\hat{Z}}^{-1}$ is bounded.\\

\begin{lm}
\label{cb}
If $S$ is T-bounded and $\det (S) \geq c >0$, then $S^{-1}$ is T-bounded.
\end{lm}
{\bf Proof:} $S$ being bounded and $\det (S) \geq c >0$ imply that $S^{-1}$ is bounded. Since $\det (\bar{\hat{Z}}S\hat{Z}^{-1}) = \det (S)$, we see that $(\bar{\hat{Z}}S\hat{Z}^{-1})^{-1} = \hat{Z}S^{-1}\bar{\hat{Z}}^{-1}$ is bounded. Consequently, $S^{-1}$ is T-bounded.
\hfill\rule{2.1mm}{2.1mm}\\

{\bf Remark:} T-boundedness can also be more generally defined for non-Hermitian matrix $G$ by requiring both $ZGZ^{-1}$ and $\bar{Z}^{-1}G\bar{Z}$ being bounded, which implies that $G$ is bounded. Similarly, corresponding $S$ is called T-bounded if both $\bar{\hat{Z}}S\hat{Z}^{-1}$ and $\bar{\hat{Z}}^{-1}S\hat{Z}$ are bounded. The two lemmas will still hold true for such general $G$ and $S$. The proofs are almost identical.\\

In the rest of this section, we will assume $G$ to be such non-Hermitian matrix. Such $G$ is T-bounded if and only if 

\[
G = \left(O\left(\min\left(\frac{|z_i|}{|z_j|},\frac{|z_j|}{|z_i|}\right)\right)\right).
\]

$G$ is called torically bounded if

\[
G = \left(O(\delta_{ij} + |z_iz_j|)\right).
\]

The matrix of a smooth toric metric on $\mathbb{C}^n$ is of this form. (With a little abuse of notation, when used with other $(n\times n)$-matrices in this section, we will use $(a_{ij})_{I\times J}$ to denote an $(n\times n)$-matrix whose $(i,j)$-entry equals to $0$ unless $i\in I$ and $j\in J$.) $G$ is called $(I\times J)$ strongly T-bounded if $G$ is a sum of a torically bounded matrix and an $(n\times n)$-matrix of the following form:

\[
\left(O\left(\frac{|z_0|^2}{|z_i||z_j|}\right)\right)_{I\times J},\ \ {\rm where}\ |z_0| \leq \min_{1\leq i \leq n}(|z_i|).
\]

The matrix of the restriction of a smooth toric metric on $\mathbb{C}^{n+1}$ to our hypersurface and the perturbed metric on $\mathbb{C}^{n+1}$ are of such form. It is easy to see that a torically bounded or strongly T-bounded matrix is T-bounded.\\  
\begin{lm}
\label{cc}
If $G$ is $(I\times J)$ strongly T-bounded and $\det (G) \geq c >0$, then $G^{-1}$ is $(I\times J)$ strongly T-bounded.
\end{lm}
{\bf Proof:} Let $G = (g_{ij})$ and $G^{-1} = (g^{ij})$. By Crammer's rule, $g^{ij} = \frac{G_{ji}}{\det(G)}$, where $G_{ji}$ is the $(ji)$-minor of $G$. When $i=j$, there is nothing to be proved. When $i\not=j$, since $G$ is bounded, through the formula of determinant, $G_{ji}$ can be expressed as sum of bounded multiples of terms like $\displaystyle \prod_{k=1}^{m-1}g_{i_ki_{k+1}}$, where $i=i_1$, $j=i_m$ and $\{i_1,\cdots,i_m\}$ are all different.\\

If $g_{i_1i_2} = O(|z_{i_1}||z_{i_2}|)$ and $\displaystyle g_{i_2i_3} = O\left(\frac{|z_0|^2}{|z_{i_2}||z_{i_3}|}\right)$, then $\displaystyle g_{i_1i_2}g_{i_2i_3} = O\left(\frac{|z_{i_1}||z_0|^2}{|z_{i_3}|}\right) = O(|z_{i_1}||z_{i_3}|)$. By induction, one can show that if $g_{i_ki_{k+1}} = O(|z_{i_k}||z_{i_{k+1}}|)$ for some $k$, then $\displaystyle \prod_{k=1}^{m-1}g_{i_ki_{k+1}} = O(|z_{i_1}||z_{i_m}|) = O(|z_i||z_j|)$.\\

On the other hand, if $\displaystyle g_{i_ki_{k+1}} = O\left(\frac{|z_0|^2}{|z_{i_k}||z_{i_{k+1}}|}\right)$ for all $k$, then $i\in I$, $j\in J$ and

\[
\prod_{k=1}^{m-1}g_{i_ki_{k+1}} \leq C|g_{i_1i_2}g_{i_{m-1}i_m}| = O\left(\frac{|z_0|^2}{|z_{i_1}||z_{i_m}|}\right) = O\left(\frac{|z_0|^2}{|z_i||z_j|}\right).
\]
\hfill\rule{2.1mm}{2.1mm}\\

Similarly, we may extend the concepts of torically bounded and strongly T-bounded to $S = \hat{G} + \hat{H}$. Let the index set $\tilde{J}$ be the union of $J$ and $n+J$ (the $n$-shift of $J$) and

\[
\tilde{S} = S\left(\begin{array}{cc}0&-I\\I&0\end{array}\right) = \left(\begin{array}{cc}G&-H\\-\bar{H}&G^t\end{array}\right).
\]

$S$ is called torically bounded or $(I,J)$ strongly T-bounded if $\tilde{S}$ is torically bounded or $(\tilde{I},\tilde{J})$ strongly T-bounded as a $(2n \times 2n)$-matrix in the previous sense, where $z_{n+i} = z_i$. (In particular, $\hat{G}$ is torically bounded or $(I,J)$ strongly T-bounded if and only if $G$ is torically bounded or $(I,J)$ strongly T-bounded in the previous sense.) Consequently we have

\begin{lm}
\label{cd}
If $S$ is $(I\times J)$ strongly T-bounded and $\det (S) \geq c >0$, then $S^{-1}$ is $(I\times J)$ strongly T-bounded.
\end{lm}
\hfill\rule{2.1mm}{2.1mm}\\

Natural examples of $G$ or $S = \hat{G} + \hat{H}$ are the matrices corresponding to $(1,1)$-forms or symplectic 2-forms. Such $(1,1)$-form or symplectic 2-form is called torically bounded, (strongly) T-bounded if the corresponding $G$ or $S = \hat{G} + \hat{H}$ is torically bounded, (strongly) T-bounded.\\ 

{\bf Remark:} Although lemmas \ref{cc} and \ref{cd} are much more refined and precise than lemmas \ref{ca} and \ref{cb}, for most application in this work, the much simpler lemmas \ref{ca} and \ref{cb} are usually sufficient.\\

\se{Estimate for the symplectic neighborhood}
For our discussion, it is crucial to establish suitable Darboux coordinates (satisfying certain estimate) for the symplectic neighborhood of the special Lagrangian local models discussed in section 2. Such coordinates will be established in this section.\\

According to proposition \ref{bd}, we may write

\[
\check{\omega}_t = (\omega - \partial\bar{\partial}\mu)|_{Y_t} = \hat{\omega}_t - id(\lambda_0|z_0|^2 - \eta)\wedge d{\rm Arg}(tp),
\]

where

\[
\hat{\omega}_t = i\sum_{k=1}^n dx_k\wedge dy_k,
\]

with the coordinate $x_k = \theta_k$ for $1\leq k \leq n$ and

\[
y_k = \lambda_k|z_k|^2 - \lambda_0|z_0|^2 - c_k,\ \ {\rm for}\ 1\leq k \leq l;
\]
\[
y_j = \rho_{j} - (\lambda_0|z_0|^2 - \eta) + \sum_{k\in I''}\frac{\lambda_{k,j}}{\lambda_k}(\lambda_k|z_k|^2 - \eta) - C^0_j,\ \ {\rm for}\ l+1\leq j \leq n.
\]

$C^0_j$ is chosen to ensure that $y_j|_{S_{r,c} \cap Y_t} = 0$. We will start with some estimates of the coordinate transformation of the coordinates $(x,y)$ and $z=(z_1,\cdots,z_n)$ of $Y_t$ near $S_{r,c} \cap Y_t$. Notice that near $S_{r,c} \cap Y_t$, $|z_k|^2 \sim (\eta+c_k)$ for $k\in I''$.\\
\begin{lm}
\label{df}
\[
\frac{\partial y_j}{\partial \log |z_k|^2} = O\left(\delta_{jk}|z_j||z_k| + |z_j|^2|z_k|^2 + |z_0|^2\right) = O\left(\min\left(|z_j|^2,|z_k|^2\right)\right),
\]
\[
\frac{\partial x}{\partial \log |z|^2}=0,\ \ \frac{\partial x}{\partial \theta}=I,
\]
\[
\frac{\partial y}{\partial \theta} = \left(\begin{array}{cc}\displaystyle\left(0\right)_{l\times l}&\displaystyle\left(O\left(|z_0|^2|z_k|\right)\right)_{l\times I'}\\\\\displaystyle\left(0\right)_{I'\times l}&\displaystyle\left(O\left(|z_0|^2|z_k|\right)\right)_{I'\times I'}\end{array}\right).
\]
\[
\frac{\partial \log |z_j|^2}{\partial y_k} = O\left(\frac{\delta_{jk}}{|z_j||z_k|} + 1 + \frac{|z_0|^2}{|z_j|^2|z_k|^2}\right) = O\left(\min\left(\frac{1}{|z_j|^2},\frac{1}{|z_k|^2}\right)\right),
\]
\[
\frac{\partial \theta}{\partial y}=0,\ \ \frac{\partial \theta}{\partial x}=I,
\]
\[
\frac{\partial \log |z|^2}{\partial x} = \left(\begin{array}{cc}\displaystyle\left(0\right)_{l\times l}&\displaystyle\left(O\left(\frac{|z_0|^2}{|z_j|^2} |z_k|\right)\right)_{l\times I'}\\\\\displaystyle\left(0\right)_{I'\times l}&\displaystyle\left(O\left(\frac{|z_0|^2}{|z_j|^2} |z_k|\right)\right)_{I'\times I'}\end{array}\right).
\]
Further more, the multi-derivatives of each non-constant term with respect to $\{\log z_k\}_{k=1}^n$ will hold same bound.\\
\end{lm}
{\bf Proof:} Recall from (\ref{be})

\[
\frac{\partial \eta}{\partial \log \kappa} = \zeta\eta,\ \ {\rm where}\ \zeta = \left(\sum_{k\in I''}\eta(c_k + \eta)^{-1}\right)^{-1}.
\]
\[
\begin{array}{c}\displaystyle\frac{\partial \eta}{\partial \theta_j} = \eta\zeta 2{\rm Re}(ip_j/p),\\\\\displaystyle\frac{\partial |z_0|^2}{\partial \theta_j} = |z_0|^2 2{\rm Re}(ip_j/p),\end{array}\ \ \ {\rm where}\ p_j = z_j\frac{\partial p}{\partial z_j}\ \ {\rm for}\ j\in I'.
\]

In the following matrix computations, we will always use $j$ as row index and $k$ as column index. It is straightforward to get

\[
\left(\begin{array}{cc}\displaystyle\frac{\partial y}{\partial \log |z|^2}&\displaystyle\frac{\partial y}{\partial \theta}\\\\\displaystyle\frac{\partial x}{\partial \log |z|^2}&\displaystyle\frac{\partial x}{\partial \theta}\end{array}\right) = \left(\begin{array}{cc}\displaystyle\frac{\partial y}{\partial \log |z|^2}&\displaystyle\frac{\partial y}{\partial \theta}\\\\0&I\end{array}\right).
\]

The inverse matrix

\[
\left(\begin{array}{cc}\displaystyle\frac{\partial \log |z|^2}{\partial y}&\displaystyle\frac{\partial \log |z|^2}{\partial x}\\\\\displaystyle\frac{\partial \theta}{\partial y}&\displaystyle\frac{\partial \theta}{\partial x}\end{array}\right) = \left(\begin{array}{cc}\displaystyle\left(\frac{\partial y}{\partial \log |z|^2}\right)^{-1}&\displaystyle-\left(\frac{\partial y}{\partial \log |z|^2}\right)^{-1}\frac{\partial y}{\partial \theta}\\\\0&I\end{array}\right).
\]

\[
\frac{\partial y}{\partial \log |z|^2} = \left(\begin{array}{cc}\displaystyle\left(\frac{\partial y_j}{\partial \log |z_k|^2}\right)_{l\times l}&\displaystyle\left(\frac{\partial y_j}{\partial \log |z_k|^2}\right)_{l\times I'}\\\\\displaystyle\left(\frac{\partial y_j}{\partial \log |z_k|^2}\right)_{I'\times l}&\displaystyle\left(\frac{\partial y_j}{\partial \log |z_k|^2}\right)_{I'\times I'}\end{array}\right)
\]
\[
= \left(\begin{array}{cc}\displaystyle\left(\delta_{jk}\lambda_k|z_k|^2 + \lambda_0|z_0|^2\right)_{l\times l}&\displaystyle\left(O(|z_j|^2|z_k|^2,|z_0|^2)\right)_{l\times I'}\\\\\displaystyle\left(\lambda_{k,j}|z_k|^2 + (\lambda_0 - \lambda_{0,j})|z_0|^2\right)_{I'\times l}&\displaystyle\left(\rho_{jk} + O(|z''|^2|z_j|^2,|z_0|^2)\right)_{I'\times I'}\end{array}\right)
\]

\[
\frac{\partial y}{\partial \theta} = \left(\begin{array}{cc}\displaystyle\left(\frac{\partial y_j}{\partial \theta_k}\right)_{l\times l}&\displaystyle\left(\frac{\partial y_j}{\partial \theta_k}\right)_{l\times I'}\\\\\displaystyle\left(\frac{\partial y_j}{\partial \theta_k}\right)_{I'\times l}&\displaystyle\left(\frac{\partial y_j}{\partial \theta_k}\right)_{I'\times I'}\end{array}\right)
\]
\[
= \left(\begin{array}{cc}\displaystyle\left(0\right)_{l\times l}&\displaystyle\left(- \lambda_0|z_0|^2 2{\rm Re}(ip_k/p)\right)_{l\times I'}\\\\\displaystyle\left(0\right)_{I'\times l}&\displaystyle\left(\left((\lambda_{0,j}-\lambda_0)|z_0|^2 + \left(1- \frac{\tilde{\Lambda}_j}{\tilde{\Lambda}}\right)\eta\zeta\right) 2{\rm Re}(ip_k/p)\right)_{I'\times I'}\end{array}\right)
\]

More precisely

\[
\left(\frac{\partial y_j}{\partial \log |z_k|^2}\right)_{l\times I'} = \left(\lambda_{j,k}|z_j|^2 - \lambda_{0,k}|z_0|^2 - \lambda_0|z_0|^2(1 + {\rm Re}(p_k/p))\right)_{l\times I'}.
\]

Notice that $|z''| \leq C|z_j|$ for $j\in I'$. It is straightforward to check that 

\[
G = Z^{-1}\frac{\partial y}{\partial \log |z|^2}\bar{Z}^{-1} = \left(\begin{array}{cc}\displaystyle\left(\delta_{jk}\lambda_k + \frac{\lambda_0|z_0|^2}{z_j\bar{z}_k}\right)_{l\times l}&\displaystyle\left(0\right)_{l\times I'}\\\\\displaystyle\left(0\right)_{I'\times l}&\displaystyle\left(\frac{\rho_{jk}}{z_j\bar{z}_k}\right)_{I'\times I'}\end{array}\right)
\]
\[
+ \left(\begin{array}{cc}\displaystyle\left(0\right)_{l\times l}&\displaystyle\left(O\left(|z_j||z_k|,\frac{|z_0|^2}{|z_j||z_k|}\right)\right)_{l\times I'}\\\\\displaystyle\left(O\left(|z_j||z_k|,\frac{|z_0|^2}{|z_j||z_k|}\right)\right)_{I'\times l}&\displaystyle\left(O\left(|z''||z_k|,\frac{|z_0|^2}{|z_j||z_k|}\right)\right)_{I'\times I'}\end{array}\right)
\]

is strongly T-bounded and 

\[
\det G = \det \left(\delta_{jk}\lambda_k + \frac{\lambda_0|z_0|^2}{z_j\bar{z}_k}\right)_{l\times l} \det \left(\frac{\rho_{jk}}{z_j\bar{z}_k}\right)_{I'\times I'} + O(|z''|,\epsilon)
\]

is bounded below by a positive constant. By lemma \ref{cc}, $G^{-1}$ is also strongly T-bounded. Hence

\[
\left(\frac{\partial \log |z_j|^2}{\partial y_k}\right) = \bar{Z}^{-1}G^{-1}Z^{-1} = \left(O\left(\frac{\delta_{jk}}{|z_j||z_k|} + 1 + \frac{|z_0|^2}{|z_j|^2|z_k|^2}\right)\right).
\]

Since

\[
\frac{\partial y}{\partial \theta} = \left(\begin{array}{cc}\displaystyle\left(0\right)_{l\times l}&\displaystyle\left(O(|z_0|^2 |z_k|)\right)_{l\times I'}\\\\\displaystyle\left(0\right)_{I'\times l}&\displaystyle\left(O(|z_0|^2 |z_k|)\right)_{I'\times I'}\end{array}\right)
\]

\[
\frac{\partial \log |z|^2}{\partial x} = - \frac{\partial \log |z|^2}{\partial y}\frac{\partial y}{\partial \theta} = \left(\begin{array}{cc}\displaystyle\left(0\right)_{l\times l}&\displaystyle\left(O\left(\frac{|z_0|^2}{|z_j|^2} |z_k|\right)\right)_{l\times I'}\\\\\displaystyle\left(0\right)_{I'\times l}&\displaystyle\left(O\left(\frac{|z_0|^2}{|z_j|^2} |z_k|\right)\right)_{I'\times I'}\end{array}\right).
\]

To prove the boundedness of multi-derivatives, notice that when $b\sim 1$, the bound of

\[
\frac{\partial}{\partial \log z_k} \left(\frac{a}{b}\right) = \frac{1}{b^2}\left(b\frac{\partial a}{\partial \log z_k} - a\frac{\partial b}{\partial \log z_k}\right)
\]

can be reduced to the bound of $\displaystyle \frac{\partial a}{\partial \log z_k}$ and $\displaystyle \frac{\partial b}{\partial \log z_k}$. Also notice that

\[
\frac{\partial |z_j|^2}{\partial \log z_k} = \delta_{jk}|z_j|^2,\ \ {\rm for}\ 1\leq j,k\leq n.
\]
\[
\frac{\partial |z_0|^2}{\partial \log z_k} = -|z_0|^2,\ \ \frac{\partial \eta}{\partial \log z_k} = 0,\ \ {\rm for}\ 1\leq k\leq l.
\]
\[
\frac{\partial |z_0|^2}{\partial \log z_k} = |z_0|^2 \frac{p_k}{p},\ \ \frac{\partial \eta}{\partial \log z_k} = \zeta\eta \left(\frac{\Lambda_k}{\Lambda} + \frac{p_k}{p} + 1\right),\ \ {\rm for}\ k\in I'.
\]
\[
\frac{\partial }{\partial \log z_k}\left(\frac{\eta}{\eta + c_j}\right) = \frac{\eta}{\eta + c_j}\left(1+\frac{\eta}{\eta + c_j}\right)\zeta \left(\frac{\Lambda_k}{\Lambda} + \frac{p_k}{p} + 1\right),\ \ {\rm for}\ k\in I'.
\]

For terms in the lemma, other than toric monomials (eigenfunctions of toric action), denominators are $\sim 1$. The log-derivatives (generators of the toric action) preserve toric monomial up to constant multiples. Since $|z_0|$ and $\eta$ never appear as denominator, above computations imply that log-derivatives will preserve the bounds of the terms in the lemma as claimed.
\hfill\rule{2.1mm}{2.1mm}\\
\begin{lm}
\label{dj}
For $I = (i_1,\cdots,i_n)$ and $J = (j_1,\cdots,j_n)$ in $\mathbb{Z}_{\geq 0}^n$, we have

\[
\left|\frac{\partial \log z_j}{\partial x_k} - i\delta_{jk}\right| \leq \frac{Cr_kr_0}{r_j},
\]
\[
\left|\frac{\partial^{|I|+|J|} \log z_j}{(\partial x)^I(\partial y)^J}\right| \leq \frac{Cr_0r_I}{r^{2J}r_j} \leq \frac{Cr_I}{r^{2J}},\ \ {\rm when}\ |I|+|J|\geq 2,\ |I| \geq 1,
\]
\[
\left|\frac{\partial^{|J|} \log z_j}{(\partial y)^J}\right| \leq \frac{C}{r^{2J}}\min\left(1, \frac{r_J^2}{r_j^2}\right) \leq \frac{C}{r^{2J}},
\]

where

\[
(\partial x)^I = \prod_{k=1}^n(\partial x_k)^{i_k},\ \ r^I = \prod_{k=1}^n r_k^{i_k},\ \ r_I = \min_{i_k>0}(r_k).
\]

In particular

\[
\left|\frac{\partial^{|I|+1} \log z_j}{(\partial x)^I\partial y_k}\right| \leq \frac{C}{r_jr_k}\min_{i_k>0}(r_k).
\]
\end{lm}
{\bf Proof:} Lemma \ref{df} implies that

\[
\left|\frac{\partial \log z_j}{\partial x_k} - i\delta_{jk}\right| \leq \frac{Cr_kr_0}{r_j},\ \ \left|\frac{\partial \log z_j}{\partial y_k} \right| \leq C\min\left(\frac{1}{r_j^2},\frac{1}{r_k^2}\right) \leq \frac{C}{r_k^2},
\]

and the multi-derivatives with respect to $\log z$ will hold same bound. Consequently

\[
\left|\frac{\partial^2 \log z_j}{\partial x_l\partial y_k} \right| = \left|\frac{\partial \log z_i}{\partial y_k}\frac{\partial}{\partial \log z_i}\left(\frac{\partial \log z_j}{\partial x_l}\right) \right| \leq C\min\left(\frac{r_l}{r_jr_k},\frac{r_l}{r_k^2}\right).
\]

\[
\left|\frac{\partial^2 \log z_j}{\partial x_k\partial x_l} \right| = \left|\frac{\partial \log z_i}{\partial x_k}\frac{\partial}{\partial \log z_i}\left(\frac{\partial \log z_j}{\partial x_l}\right) \right| \leq Cr_l.
\]

The estimates for more general $I$ and $J$ are straightforward analogue of these two estimates.
\hfill\rule{2.1mm}{2.1mm}\\

Let $\hat{\omega}_{t,s} = \hat{\omega}_t - isd(\lambda_0|z_0|^2 - \eta)\wedge d{\rm Arg}(tp)$. Then $\hat{\omega}_{t,0} = \hat{\omega}_t$ and $\hat{\omega}_{t,1} = \check{\omega}_t$.\\
\begin{lm}
\label{dd}
The matrix representing $\hat{\omega}_{t,s}$ and its inverse are strongly T-bounded.
\end{lm} 
{\bf Proof:} Since

\[
d|z_0|^2 = \sum_{k\in I'}\frac{2|z_0|^2}{|z_k|^2|p|^2}{\rm Re}(\bar{p}p_k\bar{z}_kdz_k) - \sum_{k=1}^n \frac{|z_0|^2}{|z_k|^2}d|z_k|^2 = O(\{\bar{z}_kdz_k, z_kd\bar{z}_k\}_{k=1}^n),
\]

the only term in $\omega|_{Y_t}$ that is not torically bounded is 

\[
\partial \bar{\partial} |z_0|^2 = |z_0|^2\left(\sum_{k=1}^n \frac{dz_k}{z_k} - \frac{dp}{p}\right)\left(\sum_{k=1}^n \frac{d\bar{z}_k}{\bar{z}_k} - \frac{d\bar{p}}{\bar{p}}\right).
\]

It is easy to verify that this term is strongly T-bounded. Recall that

\[
\partial\bar{\partial}\mu = \eta\partial\bar{\partial}\log \Lambda + \frac{d\eta}{d\log \kappa}\partial\log \kappa\bar{\partial}\log \kappa
\]
\[
\partial\bar{\partial}\mu|_{Y_t} = \eta\partial\bar{\partial}\log \Lambda + \zeta\eta(\partial\log \tilde{\Lambda} + d\log p)(\bar{\partial}\log \tilde{\Lambda} - d\log \bar{p}).
\]

The first term on the right is torically bounded and the second term is strongly T-bounded.

\[
\check{\omega}_t = (\omega - \partial\bar{\partial}\mu)|_{Y_t} = \sum_{k=1}^l \lambda_k\partial\bar{\partial}|z_k|^2 + \lambda_0|z_0|^2\sum_{j,k=1}^l \frac{dz_jd\bar{z}_k}{z_j\bar{z}_k} + \partial\bar{\partial}\rho + O(|z''|,\epsilon).
\]

Let $G$ be the metric matrix of $\check{\omega}_t$. We see that $G$ is strongly T-bounded and

\[
\det G = \det \left(\delta_{jk}\lambda_k + \frac{\lambda_0|z_0|^2}{z_j\bar{z}_k}\right)_{l\times l} \det \left(\frac{\rho_{jk}}{z_j\bar{z}_k}\right)_{I'\times I'} + O(|z''|,\epsilon)
\]

is bounded below by a positive constant. By lemma \ref{cc}, $G^{-1}$ is also strongly T-bounded.\\

$\hat{\omega}_{t,s}$ is non-Hermitian symplectic form and can be represented by $S = \hat{G} + \hat{H}$ (notation from section 2). It is easy to check that the perturbation $\hat{H}$ is of order $O\left(\frac{|z_0|^2}{|z_j|}\right)$ and is strongly T-bounded. Therefore $S$ and $S^{-1}$ are strongly T-bounded (lemma \ref{cd}).
\hfill\rule{2.1mm}{2.1mm}\\
\[
\frac{d\hat{\omega}_{t,s}}{ds} = d\alpha,\ \alpha = - i(\lambda_0|z_0|^2 - \eta) d{\rm Arg}(tp).
\]

Applying lemma \ref{dd}, we have

\[
H_\alpha = \imath(\alpha)\hat{\omega}_{t,s}^{-1} = (\lambda_0|z_0|^2 - \eta)\left[O\left(z''\frac{\partial}{\partial z''},\frac{\partial}{\partial z'}\right) + \sum_{k=1}^l O\left(\epsilon\frac{|z_0|}{|z_k|}\frac{\partial}{\partial z_k}\right)\right].
\]

Let $\varphi_s$ denote the flow of $-H_\alpha$. Then

\[
\frac{d\hat{\omega}_{t,s}\circ \varphi_s}{ds} = d\alpha\circ \varphi_s - d\imath(H_\alpha\circ \varphi_s)(\hat{\omega}_{t,s}\circ \varphi_s) = 0.
\]

Namely, $\varphi_s: (Y_t,\hat{\omega}_t) \rightarrow (Y_t,\hat{\omega}_{t,s})$ is a symplectomorphism. (From now on, we will use $z^0$ to denote the coordinate on $Y_t$ that was denoted by $z$.) Let $z = z^0 \circ \varphi_s$ (in another word, $z = \varphi_s(z^0)$) we have $\displaystyle\frac{dz}{ds} = -H_\alpha(z)$.\\
\begin{lm}
\label{do}
$\varphi_s|_{S_{r,c}\cap Y_t} = {\rm id}$.
\end{lm}
{\bf Proof:} Since $(\lambda_0|z_0|^2 - \eta)|_{S_{r,c}\cap Y_t}=0$, $H_\alpha|_{S_{r,c}\cap Y_t}=0$. By the uniqueness of solution of ODE, we have $\varphi_s|_{S_{r,c}\cap Y_t} = {\rm id}$.
\hfill\rule{2.1mm}{2.1mm}\\
\begin{eqnarray*}
\frac{d\log z_j}{ds} &=& -H_\alpha(\log z_j) = (\lambda_0|z_0|^2 - \eta)O\left(1+\epsilon\frac{|z_0|}{|z_j|^2}\right)\ \ {\rm for}\ 1\leq j \leq l,\\
\frac{d z_j}{ds} &=& -H_\alpha(z_j) = (\lambda_0|z_0|^2 - \eta)O(1)\ \ {\rm for}\ l+1\leq j \leq n.
\end{eqnarray*}

Hence

\[
\frac{d\log z_0}{ds} = \frac{d\log (tp)}{ds} - \sum_{j=1}^n \frac{d\log z_j}{ds} = (\lambda_0|z_0|^2 - \eta)O\left(\frac{\epsilon}{|z_0|}\right).
\]

Using the notation

\[
O(A) = AO(1),\ \ {\rm where}\ A = {\rm Diag}\left(\displaystyle\left\{\epsilon\frac{|z_0|^3}{|z_j|^2}+|z_0|^2\right\}_{j=1}^l, \displaystyle\left\{\frac{|z_0|^2}{|z_j|}\right\}_{j=l+1}^n\right)
\]

and $O(1)$ is usually a matrix, we can write

\begin{equation}
\label{da}
\frac{d\log z}{ds} = -H_\alpha(\log z) = O(A).
\end{equation}

It is easy to check that the multi-derivatives of $-H_\alpha(\log z)$ with respect to $\log z$ will hold same bound $O(A)$. Since $H_\alpha$ depends on $c$, the solution of (\ref{da}) is a function of $(s,c,z^0)$. The following lemma deals with the dependence on $c$.\\
\begin{lm}
\label{dn}
For $z = \varphi_s(z_0)$,
\[
\frac{\partial \log z}{\partial c_k} = \frac{1}{c_k + \eta}O(A),\ \ \frac{\partial \log z_0}{\partial c_k} = \frac{1}{c_k + \eta}O\left(\epsilon|z_0|\right).
\]
\end{lm}
{\bf Proof:} We may rewrite

\[
\hat{\omega}_{t,s} = \omega|_{Y_t} - \partial \bar{\partial} \mu + i(1-s)d(\lambda_0|z_0|^2 - \eta)\wedge d{\rm Arg}(tp).
\]

Applying formula (\ref{bg}), we have

\[
\frac{\partial \alpha}{\partial c_k} = - i\frac{\zeta\eta}{c_k + \eta}d{\rm Arg}(tp).
\]

\[
\frac{\partial \hat{\omega}_{t,s}}{\partial c_k} = \partial \bar{\partial} \log(c_k + \eta) + i(1-s)\frac{\zeta\eta}{c_k + \eta}d{\rm Arg}(tp)
\]
\[
= -isd\left(\frac{\zeta\eta}{c_k + \eta}\right)d{\rm Arg}(tp) -i\sum_{j\in I'} d\left(\frac{\zeta\eta}{c_k + \eta}((\log \Lambda)_j -1)\right)d\theta_j.
\]
\[
= \frac{\eta}{c_k + \eta}O\left(\frac{1}{|z_i||z_j|}\right)_{I'\times I'}.
\]

(The last term uses the matrix notation from section 2 for the 2-form.) Applying lemma \ref{dd}, we have

\[
\frac{\partial H_\alpha}{\partial c_k} = \imath\left(\frac{\partial \alpha}{\partial c_k}\right)\hat{\omega}_{t,s}^{-1} - \imath(\alpha)\hat{\omega}_{t,s}^{-1}\frac{\partial \hat{\omega}_{t,s}}{\partial c_k}\hat{\omega}_{t,s}^{-1} 
\]
\[
= \frac{\eta}{c_k + \eta}\left[O\left(z''\frac{\partial}{\partial z''},\frac{\partial}{\partial z'}\right) + \sum_{k=1}^l O\left(\epsilon\frac{|z_0|}{|z_k|}\frac{\partial}{\partial z_k}\right)\right].
\]

Consequently

\[
\frac{d}{ds}\left(\frac{\partial \log z}{\partial c_k}\right) = -\frac{\partial H_\alpha}{\partial c_k}(\log z) - \frac{\partial H_\alpha(\log z)}{\partial \log z}\frac{\partial \log z}{\partial c_k} = O(A)\left(\frac{\partial \log z}{\partial c_k}\right) + \frac{1}{c_k + \eta}O(A).
\]

Standard linear ODE estimate yields

\[
\frac{\partial \log z}{\partial c_k} = \frac{1}{c_k + \eta}O(A).
\]
\hfill\rule{2.1mm}{2.1mm}
\begin{lm}
\label{dc}
For $z = \varphi_s(z_0)$,
\[
\left(\frac{\partial \log z_j}{\partial \log z^0_k}\right) = I + O(A),\ \ \left(\frac{\partial \log z^0_j}{\partial \log z_k}\right) = I + O(A).
\]
\[
\left(\frac{\partial \log z_j}{\partial \log \bar{z}^0_k}\right) = O(A),\ \ \left(\frac{\partial \log \bar{z}^0_j}{\partial \log z_k}\right) = O(A).
\]
Further more, the multi-derivatives of each non-constant term with respect to $\{\log z^0_k\}_{k=1}^n$ ($\{\log z_k\}_{k=1}^n$) and their complex conjugates will hold same bound.
\end{lm}
{\bf Proof:} Using similar argument as in lemma \ref{df}, we see that the multi-derivatives of right hand side of equation (\ref{da}) with respect to $\{\log z_k\}_{k=1}^n$ hold the same bound. Take derivative of (\ref{da}) with respect to $\log z^0$, we have

\[
\frac{d}{ds}\left(\frac{\partial \log z}{\partial \log z^0}\right) = O(A)\left(\frac{\partial \log z}{\partial \log z^0}\right).
\]

By standard linear ODE estimate, we have

\[
\left(\frac{\partial \log z_j}{\partial \log z^0_k}\right) = I + O(A).
\]

Using induction and similar ODE estimates, we see that all the multi-derivatives of $\displaystyle \left(\frac{\partial \log z_j}{\partial \log z^0_k}\right)$ with respect to $\{\log z^0_k\}_{k=1}^n$ and their complex conjugates will hold same bound. Since

\[
\det \left(\frac{\partial \log z_j}{\partial \log z^0_k}\right) \sim 1,
\]

by standard linear algebra, we have

\[
\left(\frac{\partial \log z^0_j}{\partial \log z_k}\right) = I + O(A).
\]

Consequently

\[
\frac{\partial}{\partial \log z_i}\left(\frac{\partial \log z^0_j}{\partial \log z_k}\right) = -\left(\frac{\partial \log z^0_j}{\partial \log z_l}\right)\left(\frac{\partial^2 \log z_l}{\partial \log z^0_m\partial \log z^0_q}\frac{\partial \log z^0_q}{\partial \log z_i}\right)\left(\frac{\partial \log z^0_m}{\partial \log z_k}\right) = O(A).
\]

Similarly, by induction, all the multi-derivatives of $\displaystyle \left(\frac{\partial \log z^0_j}{\partial \log z_k}\right)$ with respect to $\{\log z_k\}_{k=1}^n$ and their complex conjugates will hold same bound. The estimates for $\displaystyle\left(\frac{\partial \log z_j}{\partial \log \bar{z}^0_k}\right)$ and $\displaystyle\left(\frac{\partial \log \bar{z}^0_j}{\partial \log z_k}\right)$ are similar. 
\hfill\rule{2.1mm}{2.1mm}\\

\se{Deformation estimate}
In this section, we will establish the estimates for the deformation from the local model discussed in section 2 to our actual Calabi-Yau hypersurface. Estimates in this section are somewhat parallel to the estimates concerning the symplectic flow $\varphi_s$ in section 4.\\ 

Consider $\mathbb{C}^{n+1}$ with coordinate $\tilde{z} = (\tilde{z}'',z')$ and toric metric with \k potential $\rho(\tilde{z}'',z')$. This can be considered as a local affine neighborhood of the toric variety. Assume the family of hypersurfaces can be locally defined as

\[
X_t = \{z_0\cdots z_n = tp(\tilde{z}'',z')\}
\]

with the \k form $\omega_t = \partial \bar{\partial} \rho|_{X_t}$. Through our discussion in this section, $t$ will be fixed. To connect with our local model, we introduce

\[
X_{t,s} = \{z_0\cdots z_n = tp(s\tilde{z}'',z')\}.
\]

$Y_t = X_{t,0}$ is our local model and $X_{t,1} = X_t$ is the actual hypersurface. Correspondingly, we may define on $\mathbb{C}^{n+1}$ the family of \k forms

\[
\tilde{\omega}_s = \partial \bar{\partial} \tilde{\rho}_s,\ \ \tilde{\rho}_s(\tilde{z}'',z') = \rho(0,z') + |s|^{-2}(\rho(s\tilde{z}'',z') - \rho(0,z')) - (1-s)\mu.
\]

Let $\tilde{\omega}_{t,s} = \tilde{\omega}_s|_{X_{t,s}}$. Clearly $\tilde{\omega}_{t,0} = \check{\omega}_t$ and $\tilde{\omega}_{t,1} = \omega_t$. $\rho$ can be locally expanded as

\[
\rho(\tilde{z}'',z') = \rho(0,z') + \sum_{k\in I''} |z_k|^2\rho^k (z') + \sum_{|I|\geq 2} |(\tilde{z}'')^I|^2\rho^I (z').
\]

Then $\tilde{\rho}_1(\tilde{z}'',z') = \rho(\tilde{z}'',z')$, $\tilde{\rho}_0(\tilde{z}'',z') = \rho(z') + |\tilde{z}''|_\lambda^2 - \mu$, where $\rho(z') = \rho(0,z')$, $\displaystyle |\tilde{z}''|_\lambda^2 = \sum_{k\in I''} \lambda_k(z') |z_k|^2$, $\lambda_k(z') = \rho^k (z')$. $\mu$ is defined in term of $\kappa(z') = \tilde{\Lambda}(z')|tp(0,z')|^2$, where $t$ is fixed. ($\mu$ here when restricted to the local model $X_{t,0}$ will coincide with the previous definition if we identify $p(z') = p(0,z')$.) We have

\[
\frac{\partial \tilde{\rho}_s}{\partial s} = v_s + \mu,\ \ {\rm where}\ v_s = \sum_{|I|\geq 2} s^{2|I|-3}|(\tilde{z}'')^I|^2\rho^I (z').
\]

Since the background metric $\tilde{\omega}_s$ on $\mathbb{C}^{n+1}$ is not fixed and depends on $s$, the Hamiltonian-gradient vector field for $(X_{t,s}, \tilde{\omega}_{t,s})$ is awkward to compute directly. Instead, we will use the Hamiltonian-gradient vector field $V$ with respect to fixed background $\tilde{\omega}_s$ (by freezing $s$), in combination with the symplectic deformation determined by varying background \k forms similar to the approach used in section 3. More precisely, $\frac{\partial \tilde{\omega}_s}{\partial s} = d \alpha_s$ where

\[
i\alpha_s = {\rm Im}(\partial \frac{\partial \tilde{\rho}_s}{\partial s}) = \sum_{k\in I''} |z_k|^2 \frac{\partial v_s}{\partial |z_k|^2}d\theta_k + \sum_{j\in I'} |z_j|^2 \frac{\partial v_s}{\partial |z_j|^2}d\theta_j + \eta {\rm Im}(\partial \log \kappa).
\]

Let $\phi_s$ be the flow determined by $V-H_{\alpha_s}$, where $H_{\alpha_s}$ is the vector field along $X_{t,s}$ satisfying

\[
\imath(H_{\alpha_s})\tilde{\omega}_{t,s} = \alpha_s|_{X_{t,s}}.
\]

\begin{lm}
\label{eb}
\[
\phi_s: (X_{t,0}, \tilde{\omega}_{t,0}) \rightarrow (X_{t,s}, \tilde{\omega}_{t,s})
\]

is a symplectomorphism.
\end{lm}
{\bf Proof:} 

\[
\frac{\partial \phi_s^*\tilde{\omega}_s}{\partial s}|_{X_{t,0}} = \phi_s^*d\imath(V-H_{\alpha_s})\tilde{\omega}_s|_{X_{t,0}} + \phi_s^*\frac{\partial\tilde{\omega}_s}{\partial s}|_{X_{t,0}}
\]
\[
= \phi_s^*d\imath(V)\tilde{\omega}_s|_{X_{t,0}} - \phi_s^*d\alpha_s|_{X_{t,0}} + \phi_s^*d\alpha_s|_{X_{t,0}} =0.
\]
\hfill\rule{2.1mm}{2.1mm}\\

{\bf Remark:} In fact, it is straightforward to check that $V-H_{\alpha_s}$ is actually the Hamiltonian-gradient vector field for the family $\{X_{t,s}, \tilde{\omega}_{t,s}\}_{s\in [0,1]}$. We will not need this fact here. We will only need lemma \ref{eb}.\\
\begin{lm}
\label{ea}
When restricted to the complex hypersurface $X_{t,s}$, where $t,s$ are constants, we have

\[
V = -\frac{2{\rm Re}(\overline{\frac{\partial tp(s\tilde{z}'',z')}{\partial s}}\nabla (z_0\cdots z_n - tp(s\tilde{z}'',z')))}{|d(z_0\cdots z_n - tp(s\tilde{z}'',z'))|^2},
\]
\[
V_t = -\frac{2{\rm Re}(\overline{p(s\tilde{z}'',z')}\nabla (z_0\cdots z_n - tp(s\tilde{z}'',z')))}{|d(z_0\cdots z_n - tp(s\tilde{z}'',z'))|^2},
\]

where $V_t$ is the Hamiltonian-gradient vector field with respect to the parameter $t$. In particular, when $X_{t,s}$ are all smooth, $V$ and $V_t$ will be smooth vector fields.
\end{lm}
{\bf Proof:} To make the notation more clear, we assume $s$ in the lemma is freezed at $s_0$. Recall that $V$ can be characterized as orthogonal to $X_{t,s_0}$ and $V(s)|_{s=s_0}=1$.  $V$ is clearly orthogonal to $X_{t,s_0}$ with respect to $\tilde{\omega}_{t,s_0}$. Since

\[
z_0\cdots z_n - tp(s\tilde{z}'',z') = 0.
\]

We have

\[
V(z_0\cdots z_n - tp(s_0\tilde{z}'',z')) = -\left.\frac{\partial tp(s\tilde{z}'',z')}{\partial s}\right|_{s=s_0}V(s)|_{s=s_0}.
\]

It is straightforward to compute

\[
V(z_0\cdots z_n - tp(s_0\tilde{z}'',z')) =-\left.\frac{\partial tp(s\tilde{z}'',z')}{\partial s}\right|_{s=s_0}.
\]

Therefore $V(s)|_{s=s_0}=1$. The derivation of the formula for $V_t$ is almost identical.
\hfill\rule{2.1mm}{2.1mm}\\

On $X_{t,s}$,

\[
d(z_0\cdots z_n - tp(s\tilde{z}'',z')) = z_0\cdots z_l\left(\sum_{k=0}^n q_k\frac{dz_k}{z_k}\right)
\]

where

\[
q_k = \left\{\begin{array}{ll}1-s\check{p}_k(s\tilde{z}'',z'),&{\rm for}\ k\in I''\\\\1-\check{p}_k(s\tilde{z}'',z'),&{\rm for}\ k\in I'\end{array}\right.,\ \ \ \ \check{p}_k = z_k\frac{\partial \check{p}}{\partial z_k},\ \check{p} = \log p.
\]

Substitute into lemma \ref{ea}, we have

\[
V = -2{\rm Re}\left(\left(\sum_{k\in I''} \check{p}_k(s\tilde{z}'',z')\right)\left(\sum_{k=0}^n \gamma_kz_k\frac{\partial}{\partial z_k}\right)\right),
\]
\[
V_t = -2{\rm Re}\left(\frac{1}{t}\left(\sum_{k=0}^n \gamma_kz_k\frac{\partial}{\partial z_k}\right)\right),
\]

where

\[
\gamma_k = \left(\sum_{j=0}^n\tilde{g}^{k\bar{j}}\frac{\bar{q}_j}{z_kz_{\bar{j}}}\right)\left(\sum_{i,j=0}^n\tilde{g}^{i\bar{j}}\frac{q_i\bar{q}_j}{z_iz_{\bar{j}}}\right)^{-1},\ \ \sum_{k=0}^n q_k \gamma_k = 1.
\]

\begin{lm}
\label{ec}
Let $(\tilde{g}_{j\bar{k}})$ and $(g_{j\bar{k}})$ denote the metric matrices of $\tilde{\omega}_s$ and $\tilde{\omega}_{t,s}$. Then both $(\tilde{g}_{j\bar{k}})$ and its inverse $(\tilde{g}^{k\bar{j}})$ are $(I'\times I')$ strongly T-bounded, and both $(g_{j\bar{k}})$ and its inverse $(g^{k\bar{j}})$ are strongly T-bounded.
\end{lm}
{\bf Proof:} Since $\partial \bar{\partial} \mu$ only depends on $z'$ and

\[
dz_0 = -\sum_{k=1}^n\frac{z_0q_k}{z_kq_0}dz_k,\ \ {\rm on}\ X_{t,s}.
\]

It is straightforward to verify from the definitions that $(\tilde{g}_{j\bar{k}})$ is $(I'\times I')$ strongly T-bounded, and $(g_{j\bar{k}})$ is strongly T-bounded.\\

Since $\partial \bar{\partial} \mu = O(\epsilon)$, $\det(\tilde{g}_{j\bar{k}})>C>0$. $\det(g_{j\bar{k}})>C>0$ for similar reason. By lemma \ref{cc}, $(\tilde{g}^{k\bar{j}})$ is $(I'\times I')$ strongly T-bounded, and $(g^{k\bar{j}})$ is strongly T-bounded.
\hfill\rule{2.1mm}{2.1mm}\\
\begin{lm}
\label{ef}
\[
V = \sum_{k=0}^n O\left(\frac{|z''||z_0|^2}{|z_k|^2}\right)z_k\frac{\partial}{\partial z_k},
\]
\[
H_{\alpha_s} = O\left(\left\{\max\left(|z''|^2, \epsilon^2\frac{|z_0|^2}{|z_k|^2}\right)z_k\frac{\partial}{\partial z_k}\right\}_{k=1}^l\right) + O\left(\left\{\frac{|z_0|^2}{|z_k|}\frac{\partial}{\partial z_k}\right\}_{k=l+1}^n\right).
\]
\end{lm}
{\bf Proof:} Lemma \ref{ec} implies that $\displaystyle \gamma_k = O\left(\frac{|z_0|^2}{|z_k|^2}\right)$. This give us the estimate for $V$. It is straightforward to compute that

\[
\partial v_s = \sum_{k=1}^l O(|z''|^2\bar{z}_k dz_k) + \sum_{k=l+1}^n O(|z''|^4\bar{z}_k dz_k),
\]
\[
\partial\mu = \eta\partial\log \kappa = \eta(\partial\log \tilde{\Lambda} + d\log p) = \sum_{k=l+1}^nO\left(\eta \frac{dz_k}{z_k}\right).
\]

Applying lemma \ref{ec} again, we get the estimate for $H_{\alpha_s}$.
\hfill\rule{2.1mm}{2.1mm}\\

Let $z^0 = z|_{X_{t,0}}$ and $z^s = z\circ\phi_s|_{X_{t,0}}$ (in another word, $z^s = \phi_s(z^0)$) be functions on $X_{t,0}$ (for simplicity of notation, we will use $z$ to denote $z^s$) and

\[
B = {\rm Diag}\left(\left\{\max\left(|z''|^2, (\epsilon^2 + |z''|)\frac{|z_0|^2}{|z_j|^2}\right)\right\}_{j=1}^l, \left\{\frac{|z_0|^2}{|z_j|^2}\right\}_{j=l+1}^n\right).
\]

Applying lemma \ref{ef}, we have

\begin{equation}
\label{ed}
\frac{d\log z}{ds} = O(B),\ \ \frac{d\log z_0}{ds} = O(|z''|,\epsilon^2).
\end{equation}

\begin{lm}
\label{ee}
For $z = \phi_s(z_0)$,
\[
\log z - \log z^0 = O(B),\ \ \left(\frac{\partial \log z_j}{\partial \log z^0_k}\right) = I + O(B),\ \ \left(\frac{\partial \log z^0_j}{\partial \log z_k}\right) = I + O(B).
\]
Further more, the multi-derivatives of each non-constant term with respect to $\{\log z^0_k\}_{k=1}^n$ ($\{\log z_k\}_{k=1}^n$) and their complex conjugates will hold same bound.
\end{lm}
{\bf Proof:} The proof here is identical to the proof of lemma \ref{dc} by replacing $A$ with $B$.
\hfill\rule{2.1mm}{2.1mm}\\
\begin{lm}
\label{eg}
For $z = \phi_s(z_0)$,
\[
\frac{\partial \log z}{\partial c_k} = \frac{1}{c_k + \eta}O(B),\ \ \frac{\partial \log z_0}{\partial c_k} = \frac{1}{c_k + \eta}O\left(|z''|,\epsilon^2\right).
\]
\end{lm}
{\bf Proof:} Applying formula (\ref{bg}), we have

\[
\frac{\partial \alpha_s}{\partial c_k} = i{\rm Im}(\partial\log(c_k + \eta)) = i\frac{\zeta\eta}{c_k + \eta}{\rm Im}(\partial\log \kappa) = \frac{\eta}{c_k + \eta}\sum_{j=l+1}^nO\left(\frac{dz_j}{z_j}\right).
\]

As in lemma \ref{dn}

\[
\frac{\partial \tilde{\omega}_s}{\partial c_k} = (1-s)\partial \bar{\partial} \log(c_k + \eta) = \frac{\eta}{c_k + \eta}O\left(\frac{1}{|z_i||z_j|}\right)_{I'\times I'}.
\]
\[
\frac{\partial \tilde{\omega}_{t,s}}{\partial c_k} = \left.\frac{\partial \tilde{\omega}_s}{\partial c_k}\right|_{X_{t,s}} = \frac{\eta}{c_k + \eta}O\left(\frac{1}{|z_i||z_j|}\right)_{I'\times I'}.
\]

(The last term uses the matrix notation from section 2 for the $(1,1)$-form.)

\[
\frac{\partial \gamma_j}{\partial c_k} = \left\{
\begin{array}{ll}\displaystyle\frac{\epsilon^2}{c_k + \eta}\frac{|z_0|^2}{|z_j|^2},&j\in I'\\\displaystyle\frac{\epsilon^2}{c_k + \eta}|z_0|^2,&j\in I''\end{array}\right.
\]
\[
\frac{\partial V}{\partial c_k} = \frac{\epsilon^2}{c_k + \eta}\left(\sum_{j=l+1}^n O\left(\frac{|z''||z_0|^2}{|z_j|^2}\right)z_j\frac{\partial}{\partial z_j} + \sum_{j=0}^l O\left(|z''||z_0|^2\right)z_j\frac{\partial}{\partial z_j}\right)
\]

Applying lemma \ref{ec}, we have

\[
\frac{\partial H_\alpha}{\partial c_k} = \imath\left(\frac{\partial \alpha_s}{\partial c_k}\right)\tilde{\omega}_{t,s}^{-1} - \imath(\alpha)\tilde{\omega}_{t,s}^{-1}\frac{\partial \tilde{\omega}_{t,s}}{\partial c_k}\tilde{\omega}_{t,s}^{-1} 
\]
\[
= \frac{\eta}{c_k + \eta}\left[\sum_{j=l+1}^n O\left(\frac{1}{|z_j|^2}\right)z_j\frac{\partial}{\partial z_j} + \sum_{j=1}^l O\left(\max\left(1,\frac{\epsilon^2}{|z_j|^2}\right)\right)z_j\frac{\partial}{\partial z_j}\right].
\]

Consequently, (as in lemma \ref{dn})

\[
\frac{d}{ds}\left(\frac{\partial \log z}{\partial c_k}\right) = O(B)\left(\frac{\partial \log z}{\partial c_k}\right) + \frac{1}{c_k + \eta}O(B).
\]

Standard linear ODE estimate yields

\[
\frac{\partial \log z}{\partial c_k} = \frac{1}{c_k + \eta}O(B).
\]
\hfill\rule{2.1mm}{2.1mm}\\

{\bf Remark:} It is easy to see that $B$ dominates $A$. Lemmas \ref{dc} and \ref{ee}, lemmas \ref{dn} and \ref{eg} can be combined. Consequently, the estimates in lemmas \ref{ee} and \ref{eg} are also true for $z = \psi_s(z^0)$, where $\psi_s = \phi_s \circ \varphi_1$.\\

\se{Basic estimate}
Having established the symplectic neighborhood estimate for the local model $(Y_t, \check{\omega}_t)$ in section 4 and the estimate for the coordinate under deformation in section 5, we will apply them here to work out the basic estimates necessary to apply the implicit function theorem to our deformation construction. Recall the basic setting from \cite{sl1,sl2}.\\
\[
\Omega = \prod_{i=0}^n \frac{dz_i}{z_i},\ \ \frac{dz_i}{z_i} = \frac{\partial \log z_i}{\partial y_k}dy_k + \frac{\partial \log z_i}{\partial x_k}dx_k.
\]

\[
\Omega_{t,s} = \imath(V_t)\Omega|_{X_{t,s}} = -\frac{1}{q_0}\prod_{i=1}^n \frac{dz_i}{z_i} = \eta_s\prod_{i=1}^n \left(dx_i + (U_s)_{ij}dy_j\right).
\]

\[
\eta_s = -\frac{1}{q_0}\det\left(\frac{\partial \log z_i}{\partial x_k}\right),\ \ (U_s)_{ij} = \left(\frac{\partial \log z_k}{\partial x_i}\right)^{-1}\left(\frac{\partial \log z_k}{\partial y_j}\right). 
\]

Let $z = \psi_s(z^0)$ and

\[
(U^1_s)_{ij} = \frac{\partial \log z_i}{\partial \log z^0_j},\ \ (\tilde{U}^1_s)_{ij} = \frac{\partial \log z_i}{\partial \log \bar{z}^0_j},\ \ (U^2)_{ij} = \frac{\partial \log z^0_i}{\partial x_j},\ \ (U^3)_{ij} = \frac{\partial \log z^0_i}{\partial y_j}. 
\]

Then
\begin{equation}
\label{fe}
\eta_s = -\frac{1}{q_0}\det(U^1_sU^2 +\tilde{U}^1_s\bar{U}^2),\ U_s = (U^1_sU^2 +\tilde{U}^1_s\bar{U}^2)^{-1}(U^1_sU^3 +\tilde{U}^1_s\bar{U}^3)
\end{equation}

Let $L_{t,0} = S_{r,c}\cap X_{t,0}$. Since $\tilde{g}_0|_{L_{t,0}}$ is not flat, it is convenient to introduce a background flat metric $\check{g}$ on $L_{t,0}$, which is quasi-isometric to $\tilde{g}_0|_{L_{t,0}}$, to measure the estimates. Let $\nu_i = r^0_i$ for $i\in I'$, and $\nu_i = \sqrt{\check{\eta} + c_i}$ for $i\in I''$. Define $\displaystyle \check{g} = \sum_{i=1}^n \nu_i^2dx_i^2$. Since in our normal region $r^0_i = |z^0_i| \sim \nu_i$ near $L_{t,0}$ for all $i$, $\check{g}$ is quasi-isometric to $\tilde{g}_0|_{L_{t,0}}$. From this section on, we will further require our {\bf normal region} to satisfy $|z^0_i|^{\frac{3}{2}} \leq C|z^0_1|$ for $1\leq i \leq l$ and $|z^0_i|^{\frac{3}{2}} \geq C|z^0_0|$ for $i \in I'$. Equivalently, $\nu_i^{\frac{3}{2}} \leq C\nu_1$ for $1\leq i \leq l$ and $\nu_i^{\frac{3}{2}} \geq C\nu_0$ for $i \in I'$. Under such conditions, we have $\epsilon \leq C|z^0_0|^{\frac{1}{3}} \sim \nu_0^{\frac{1}{3}}$ and $|z''| \leq C|z^0_1|^{\frac{2}{3}} \sim \nu_1^{\frac{2}{3}}$.\\
\[
B = {\rm Diag}\left(\displaystyle\left\{\nu_1^{\frac{1}{3}}\max\left(\nu_1, \frac{|z_0|}{|z_j|}\right)\right\}_{j=1}^n\right).
\]

Recall the transition estimates of norms on thin torus (lemma 3.5 of \cite{sl2}).\\
\begin{lm}
\label{fg}
\[
\left[\frac{\partial h}{\partial x_k}\right]_{C^0} \leq C\nu_k^{2+\alpha}[h]_{C^{2,\alpha}},\ \left[\frac{\partial^2 h}{\partial x_j\partial x_k}\right]_{C^0} \leq C\nu_j\nu_k\min(\nu_j^\alpha,\nu_k^\alpha)[h]_{C^{2,\alpha}}.
\]
\[
\left[\frac{\partial h}{\partial x_k}\right]_{C^\alpha} \leq C\nu_k^2[h]_{C^{2,\alpha}},\ \left[\frac{\partial^2 h}{\partial x_j\partial x_k}\right]_{C^\alpha} \leq C\nu_j\nu_k[h]_{C^{2,\alpha}}.
\]
\end{lm}
\hfill\rule{2.1mm}{2.1mm}\\

A Lagrangian torus $L$ near $L_{t,0}$ can be expressed as $(x,y) = (x,\frac{\partial h}{\partial x})$ under the Darboux coordinate from section 4. As in \cite{sl1,sl2},

\[
F(h,s) = {\rm Im}\left(\log \Omega_{t,s}|_L\right) = {\rm Im}\left(\log \eta_s\left(x, \frac{\partial h}{\partial x}\right)+ \log \det\left(I + U_s\frac{\partial^2 h}{\partial x^2}\right)\right)
\] 

defines a map $F: {\cal B}_1\times \mathbb{R} \rightarrow {\cal B}_2$, where ${\cal B}_1 = C^{2,\alpha}(L_{t,0})$ and ${\cal B}_2 = C^\alpha(L_{t,0})$ are Banach spaces. We intend to apply implicit function theorem to $F$ to construct the family of generalized special Lagrangians $L_{t,s}$ with respect to $(X_{t,s},\tilde{\omega}_{t,s},\Omega_{t,s})$.\\

In this section, we will always assume that $[h]_{C^{2,\alpha}} \leq C$ for a constant $C>0$ that will be determined later in theorem \ref{fd} when we try to apply the implicit function theorem. Then according to lemma \ref{fg}, we have $|\frac{\partial h}{\partial x_k}|_{C^0} \leq C\nu_k^{2+\alpha}$. For this reason, we will always restrict our estimates in this section to the neighborhood of $L_{t,0}$, where $|y_k| \leq C\nu_k^{2+\alpha}$.\\
\begin{lm}
\label{fc}
\[
\left[(U^1_s)_{ij} \left(x,\frac{\partial h}{\partial x}\right) - U^1_0 \left(x,0\right)\right]_{C^0} \leq C\nu_1^{\frac{1}{3}} \min\left(1,\frac{\nu_j}{\nu_i}\right).
\] 
\[
\left[(U^1_s)_{ij} \left(x,\frac{\partial h}{\partial x}\right)\right]_{C^{\alpha}} \leq C\left(\nu_1^{\frac{1}{3}-\alpha} + \nu_1^{\frac{1}{3}}[h]_{C^{2,\alpha}}\right) \min\left(1,\frac{\nu_j}{\nu_i}\right).
\] 
\[
\left[\left(\frac{\partial U^1_s}{\partial y_k}\right)_{ij} \left(x,\frac{\partial h}{\partial x}\right)\right]_{C^{\alpha}} \leq \frac{C}{\nu_k^2}\left(\nu_1^{\frac{1}{3}-\alpha} + \nu_1^{\frac{1}{3}}[h]_{C^{2,\alpha}}\right) \min\left(1,\frac{\nu_j}{\nu_i}\right).
\]

Same estimates hold for $\tilde{U}^1_s$.
\end{lm}
{\bf Proof:} These estimates are direct corollary of lemmas \ref{dc} and \ref{ee} (also see the remark after \ref{eg}). The last two estimates will also need lemma \ref{dj}. We will prove the second estimate to illustrate.

\[
\left[U^1_s \left(x,\frac{\partial h}{\partial x}\right)\right]_{C^{\alpha}} 
\] 
\[
\leq C\left[\frac{\partial U^1_s}{\partial x}\right]_{C^0} \max_{1\leq l\leq n} \left(|z^0_l|^{-\alpha}\right)
+ C\max_{1\leq k\leq n} \left(\left[\frac{\partial U^1_s}{\partial y_k}\right]_{C^0}\left[\frac{\partial h}{\partial x_k}\right]_{C^{\alpha}}\right) 
\]

Applying lemmas \ref{dj}, \ref{dc} and \ref{ee}, we have

\[
\left[\left(\frac{\partial U^1_s}{\partial x}\right)_{ij}\right]_{C^0} = \left[\frac{\partial U^1_s}{\partial \log z^0} \frac{\partial \log z^0}{\partial x} + \frac{\partial U^1_s}{\partial \log \bar{z}^0} \frac{\partial \log \bar{z}^0}{\partial x}\right]_{C^0} \leq C\nu_1^{\frac{1}{3}} \min\left(1,\frac{\nu_j}{\nu_i}\right),
\]
\[
\left[\left(\frac{\partial U^1_s}{\partial y_k}\right)_{ij}\right]_{C^0} = \left[\frac{\partial U^1_s}{\partial \log z^0} \frac{\partial \log z^0}{\partial y_k} + \frac{\partial U^1_s}{\partial \log \bar{z}^0} \frac{\partial \log \bar{z}^0}{\partial y_k}\right]_{C^0} \leq \frac{C\nu_1^{\frac{1}{3}}}{\nu_k^2} \min\left(1,\frac{\nu_j}{\nu_i}\right).
\]

Combining these estimates together with lemma \ref{fg}, we have

\[
\left[(U^1_s)_{ij} \left(x,\frac{\partial h}{\partial x}\right)\right]_{C^{\alpha}} \leq C\left(\nu_1^{\frac{1}{3}-\alpha} + \nu_1^{\frac{1}{3}}[h]_{C^{2,\alpha}}\right) \min\left(1,\frac{\nu_j}{\nu_i}\right).
\] 
\hfill\rule{2.1mm}{2.1mm}\\
In this paper, $C^{2,0}$ is the same as $C^2$.\\
\begin{lm}
\label{fh}
For $\beta = 0$ or $\alpha$
\[
\left[(U^2)_{ij} \left(x,\frac{\partial h}{\partial x}\right) - (U^2)_{ij} \left(x,0\right)\right]_{C^{\beta}} \leq C \nu_j^{1-\beta}[h]_{C^{2,\beta}}.
\] 
\[
\left[\left(\frac{\partial U^2}{\partial y_k}\right)_{ij} \left(x,\frac{\partial h}{\partial x}\right)\right]_{C^{\beta}} \leq \frac{C\nu_j^{1-\beta}}{\nu_k^2}.
\]
\end{lm}
{\bf Proof:} Lemma \ref{dj} implies that

\[
\left[\left(\frac{\partial U^2}{\partial y_k}\right)_{ij} \left(x,\frac{\partial h}{\partial x}\right)\right]_{C^0} \leq \frac{C\nu_j}{\nu_k^2}.
\]
\[
\left[\left(\frac{\partial^2 U^2}{\partial y_k\partial y_l}\right)_{ij}\right]_{C^0} \leq \frac{C\nu_j}{\nu_k^2\nu_l^2},\ \ \left[\left(\frac{\partial^2 U^2}{\partial x_l\partial y_k}\right)_{ij}\right]_{C^0} \leq \frac{C\nu_j}{\nu_k^2}.
\]

\[
\left[\left(\frac{\partial U^2}{\partial y_k}\right)_{ij} \left(x,\frac{\partial h}{\partial x}\right)\right]_{C^{\alpha}} \leq C \max_{1\leq l\leq n} \left(\left[\left(\frac{\partial^2 U^2}{\partial x_l\partial y_k}\right)_{ij}\right]_{C^0} |z^0_l|^{-\alpha}\right) 
\]
\[
+ C\max_{1\leq l\leq n} \left(\left[\left(\frac{\partial^2 U^2}{\partial y_k\partial y_l}\right)_{ij}\right]_{C^0}\left[\frac{\partial h}{\partial x_l}\right]_{C^{\alpha}}\right) \leq \frac{C}{\nu_k^2}(\nu_j^{1-\alpha} + \nu_j[h]_{C^{2,\alpha}})
\]

\[
(U^2)_{ij} \left(x,\frac{\partial h}{\partial x}\right) - (U^2)_{ij} \left(x,0\right) = \int_0^1 \left(\frac{\partial U^2}{\partial y_k}\right)_{ij} \left(x,\tau\frac{\partial h}{\partial x}\right)\frac{\partial h}{\partial x_k}d\tau.
\]

For $\beta = 0$ or $\alpha$

\[
\left[(U^2)_{ij} \left(x,\frac{\partial h}{\partial x}\right) - (U^2)_{ij} \left(x,0\right)\right]_{C^{\beta}} 
\]
\[
\leq \int_0^1 \left(\left[\left(\frac{\partial U^2}{\partial y_k}\right)_{ij} \left(x,\tau\frac{\partial h}{\partial x}\right)\right]_{C^0}\left[\frac{\partial h}{\partial x_k}\right]_{C^{\beta}} + \left[\left(\frac{\partial U^2}{\partial y_k}\right)_{ij} \left(x,\tau\frac{\partial h}{\partial x}\right)\right]_{C^{\beta}}\left[\frac{\partial h}{\partial x_k}\right]_{C^0}\right)d\tau.
\]
\[
\leq \frac{C\nu_j^{1-\beta}}{\nu_k^2}\left(\left[\frac{\partial h}{\partial x_k}\right]_{C^{\beta}} + \left[\frac{\partial h}{\partial x_k}\right]_{C^0}\right) \leq C\nu_j^{1-\beta} [h]_{C^{2,\beta}}.
\]
\hfill\rule{2.1mm}{2.1mm}
\begin{lm}
\label{fk}
For $\beta = 0$ or $\alpha$
\[
\left[(U^3)_{ij} \left(x,\frac{\partial h}{\partial x}\right) - (U^3)_{ij} \left(x,0\right)\right]_{C^{\beta}} \leq \frac{C}{\nu_i\nu_j} [h]_{C^{2,\beta}}.
\] 
\[
\left[(U^3)_{ij} \left(x,\frac{\partial h}{\partial x}\right)\right]_{C^{\beta}} \leq \frac{C}{\nu_i\nu_j}.
\] 
\[
\left[(\frac{\partial U^3}{\partial y_k})_{ij} \left(x,\frac{\partial h}{\partial x}\right)\right]_{C^{\beta}} \leq \frac{C}{\nu_i\nu_j\nu_k^2}.
\]
\end{lm}
{\bf Proof:} The proof is similar to lemma \ref{fh} by reducing all estimates to lemma \ref{dj}. The appearance of the factor $\nu_i\nu_j$ instead of $\nu_j^2$ is due to the estimate

\[
\left|\frac{\partial \log z^0_i}{\partial y_j}\right| \leq \frac{C}{\nu_i\nu_j}
\]

that is implied by lemma \ref{dj}.
\hfill\rule{2.1mm}{2.1mm}\\

Lemmas \ref{fc}, \ref{fh} and \ref{fk} together applying to (\ref{fe}) imply the following.\\
\begin{lm}
\label{fm}
For $\beta = 0$ or $\alpha$ 
\[
\left[\log \eta_s \left(x,\frac{\partial h}{\partial x}\right) - \log \eta_0 \left(x,0\right)\right]_{C^{\beta}} \leq C\left(\nu_1^{\frac{1}{3}-\beta} + [h]_{C^{2,\beta}}\right).
\] 
\[
\left[\frac{1}{\eta_s}\frac{\partial \eta_s}{\partial y_k} \left(x,\frac{\partial h}{\partial x}\right) - \frac{1}{\eta_0}\frac{\partial \eta_0}{\partial y_k} \left(x,0\right)\right]_{C^{\beta}} \leq \frac{C}{\nu_k^2}\left(\nu_1^{\frac{1}{3}-\beta} + [h]_{C^{2,\beta}}\right).
\] 
\[
\left[(U_s)_{ij} \left(x,\frac{\partial h}{\partial x}\right) - (U_0)_{ij} \left(x,0\right)\right]_{C^0} \leq \frac{C}{\nu_i\nu_j}(\nu_1^{\frac{1}{3}} + [h]_{C^2}).
\] 
\[
\left[(U_s)_{ij} \left(x,\frac{\partial h}{\partial x}\right) - (U_0)_{ij} \left(x,0\right)\right]_{C^{\alpha}} \leq \frac{C}{\nu_i\nu_j}\left(\nu_1^{\frac{1}{3}-\alpha} + [h]_{C^{2,\alpha}}\right) + \frac{C}{\nu_i^{1-\alpha}\nu_j^{1+\alpha}}[h]_{C^2}.
\] 
\[
\left[(U_s)_{ij} \left(x,\frac{\partial h}{\partial x}\right)\right]_{C^{\beta}} \leq \frac{C}{\nu_i^{1-\beta}\nu_j^{1+\beta}}.
\] 
\[
\left[\left(\frac{\partial U_s}{\partial y_k}\right)_{ij} \left(x,\frac{\partial h}{\partial x}\right)\right]_{C^{\beta}} \leq \frac{C}{\nu_i^{1-\beta}\nu_j^{1+\beta}\nu_k^2}.
\]
\end{lm}
{\bf Proof:} Only the estimates of $\eta_s$ need some additional comments. In $\eta_s$, other than those partial derivative terms that can be dealt with using lemmas \ref{fc}, \ref{fh} and \ref{fk}, there is an additional factor $q_0 = 1- s\check{p}_0$, which need the following estimate

\[
[\log q_0]_{C^\beta} \leq C(\nu_0)^{1-\beta} \leq C\nu_1^{\frac{1}{3}-\beta}
\]

that is straightforward to verify.
\hfill\rule{2.1mm}{2.1mm}\\
\begin{lm}
\label{fi}
Recall $a_s^{ij}$ and $b_s^i$ from \cite{sl1,sl2}. For $\beta =0$ or $\alpha$,
\[
[a_s^{ij} - a^{ij}]_{C^{\beta}} \leq \frac{C}{\nu_i\nu_j^{1+\beta}}\left(\nu_1^{\frac{1}{3}-\beta} + [h]_{C^{2,\beta}}\right).
\]
\[
[b_s^i - b^i]_{C^{\beta}} \leq \frac{C}{\nu_i^2}\left(\nu_1^{\frac{1}{3}-\beta} + [h]_{C^{2,\beta}}\right).
\]
\end{lm}
{\bf Proof:} Let $\Upsilon = {\rm Diag}(\nu_1,\cdots,\nu_n)$, $(a_s^{ij})^\Upsilon = \Upsilon(a_s^{ij})\Upsilon$, $U_s^\Upsilon = \Upsilon U_s\Upsilon$, $\left(\frac{\partial^2 h}{\partial x^2}\right)^\Upsilon = \Upsilon^{-1}\frac{\partial^2 h}{\partial x^2}\Upsilon^{-1}$. Then

\[
(a_s^{ij})^\Upsilon = {\rm Im}\left(\left(I + U_s^\Upsilon\left(\frac{\partial^2 h}{\partial x^2}\right)^\Upsilon\right)^{-1} U_s^\Upsilon\right)\left(x,\frac{\partial h}{\partial x}\right),
\]
\[
b_s^i(x) = {\rm Im}\left(\frac{1}{\eta_s}\frac{\partial \eta_s}{\partial y_i} + {\rm Tr}\left(\left(I + U_s^\Upsilon\left(\frac{\partial^2 h}{\partial x^2}\right)^\Upsilon\right)^{-1}\frac{\partial U_s^\Upsilon}{\partial y_i}\left(\frac{\partial^2 h}{\partial x^2}\right)^\Upsilon\right)\right)\left(x,\frac{\partial h}{\partial x}\right).
\]

Estimates in lemma \ref{fm} imply that for $\beta = 0$ or $\alpha$

\[
\left[U_s^\Upsilon \Upsilon^\beta \left(x,\frac{\partial h}{\partial x}\right)\right]_{C^{\beta}} \leq C,
\]
\[
\left[U_s^\Upsilon \Upsilon^\beta \left(x,\frac{\partial h}{\partial x}\right) - U_0^\Upsilon \Upsilon^\beta \left(x,0\right)\right]_{C^{\beta}} \leq C\left(\nu_1^{\frac{1}{3}-\beta} + [h]_{C^{2,\beta}}\right).
\]

Also notice that by lemma \ref{fg}

\[
\left[\left(\frac{\partial^2 h}{\partial x^2}\right)^\Upsilon\right]_{C^{\beta}} = [h]_{C^{2,\beta}},\ \ \left[\Upsilon^{-\alpha}\left(\frac{\partial^2 h}{\partial x^2}\right)^\Upsilon\right]_{C^0} \leq C [h]_{C^{2,\alpha}}.
\]

Consequently

\[
\left[(a_s^{ij})^\Upsilon\Upsilon^\beta - (a^{ij})^\Upsilon\Upsilon^\beta\right]_{C^{\beta}} \leq C\left(\nu_1^{\frac{1}{3}-\beta} + [h]_{C^{2,\beta}}\right). 
\]

The estimate for $b_s^i$ is similar by applying lemma \ref{fm}.
\hfill\rule{2.1mm}{2.1mm}\\
\begin{lm}
\label{fa}
\[
\left|\frac{\partial a^{ij}}{\partial x_k}\right| \leq \frac{C\nu_k}{\nu_i\nu_j},\ \ |b^i| + \left|\frac{\partial b^i}{\partial x_k}\right| \leq \frac{C\nu_0}{\nu_i^2}.
\]
Consequently
\[
[a^{ij}]_{C^\alpha} \leq \frac{C}{\nu_i\nu_j},\ \ [b^i]_{C^\alpha} \leq \frac{C(\nu_0)^{1-\alpha}}{\nu_i^2}.
\]

\end{lm}
{\bf Proof:} Recall that

\[
a^{ij} = {\rm Im}(\hat{a}^{ij}),\ (\hat{a}^{ij}(x)) = (U^2)^{-1} U^3(x,0),\ \ (U^2)_{jk} = \frac{\partial \log z^0_j}{\partial x_k},\ (U^3)_{jk} = \frac{\partial \log z^0_j}{\partial y_k},
\]
\[
b^i = {\rm Im}(\hat{b}^i),\ \hat{b}^i(x) = \frac{\partial \log \eta_0}{\partial y_i}(x,0) = {\rm Tr}\left((U^2)^{-1}\frac{\partial U^2}{\partial y_i}\right)(x,0),\ \ \eta_0 = \det U^2.
\]

\[
\frac{\partial (\hat{a}^{ij})}{\partial x_k} = -(U^2)^{-1}\frac{\partial U^2}{\partial x_k}(U^2)^{-1}U^3 + (U^2)^{-1}\frac{\partial U^3}{\partial x_k}.
\]

$\Upsilon(U^2)^{-1}\Upsilon^{-1}$ and $\Upsilon U^3 \Upsilon$ are bounded, $\displaystyle\Upsilon\frac{\partial U^2}{\partial x_k}\Upsilon^{-1}$ and $\displaystyle\Upsilon\frac{\partial U^3}{\partial x_k}\Upsilon$ are bounded by $\nu_k$. These give the first estimate.
\hfill\rule{2.1mm}{2.1mm}\\
\begin{prop}
\label{fb}
Assume that $t$ is small and the torus $(L_{t,0}, g_{t,0}|_{L_{t,0}})$ has bounded diameter. Then there exists a constant $C$ such that

\[
\|\delta h\|_{{\cal B}_1} \leq C \left\|\left(\frac{\partial F}{\partial h}(0,0)\right)\delta h\right\|_{{\cal B}_2}.
\]
\end{prop}
{\bf Proof:} One only need to show that $C$ is independent of the thinness of the torus. Recall that

\[
\delta_h F(0,0) = \frac{\partial F}{\partial h}(0,0) \delta h = a^{ij}\frac{\partial^2 \delta h}{\partial x_i\partial x_j} + b^i\frac{\partial\delta h}{\partial x_i}.
\]

We may rewrite this as

\[
a_\Upsilon^{ij}\frac{\partial^2 \delta h}{\partial x^\Upsilon_i\partial x^\Upsilon_j} = \frac{\partial F}{\partial h}(0,0) \delta h - b^i\frac{\partial\delta h}{\partial x_i},
\]

where $a_\Upsilon^{ij} = \nu_i\nu_ja^{ij}$, $x^\Upsilon_i = x_i \nu_i^{-1}$. Under the coordinate $x^\Upsilon$, $\displaystyle \check{g} = \sum_{i=1}^n (dx^\Upsilon_i)^2$.\\

According to (3.1) in \cite{sl1}, $a^{ij}(x) = (g_{t,0}|_{L_{t,0}})^{ij}$, which is clearly uniformly elliptic with respect to $\check{g}$. Lemma \ref{fa} implies that $|a_\Upsilon^{ij}|_{C^\alpha}$ is bounded.\\

It is easy to see that there exists a finite covering map $\pi: \tilde{L} \rightarrow L_{t,0}$ such that $(\tilde{L}, \tilde{g} = \pi^*\check{g})$ is of bounded geometry. We will use those symbols with ``$\ \tilde{ }\ $" to denote the corresponding pullback objects by $\pi$. By standard Schauder estimate, we have

\[
\|\tilde{\delta h}\|_{\tilde{{\cal B}}_1} \leq C \left\|\left(\frac{\partial \tilde{F}}{\partial \tilde{h}}(0,0)\right)\tilde{\delta h} - \tilde{b}^i\frac{\partial\tilde{\delta h}}{\partial \tilde{x}_i}\right\|_{\tilde{{\cal B}}_2},
\]

where the constant $C$ only depends on the dimension. By lemmas \ref{fg} and \ref{fa},

\[
\left\|\tilde{b}^i\frac{\partial\tilde{\delta h}}{\partial \tilde{x}_i}\right\|_{\tilde{{\cal B}}_2} \leq C|\tilde{b}^i|_{C^\alpha}\|\tilde{\delta h}\|_{\tilde{{\cal B}}_1} \leq C\nu_0^{1-\alpha}\|\tilde{\delta h}\|_{\tilde{{\cal B}}_1}.
\]

Since $\nu_0$ is small when $t$ is small, this term can be absorbed by $\|\tilde{\delta h}\|_{\tilde{{\cal B}}_1}$. We have

\[
\|\tilde{\delta h}\|_{\tilde{{\cal B}}_1} \leq C \left\|\left(\frac{\partial \tilde{F}}{\partial \tilde{h}}(0,0)\right)\tilde{\delta h}\right\|_{\tilde{{\cal B}}_2}.
\]

Since $\tilde{\delta h}$ and $\left(\frac{\partial \tilde{F}}{\partial \tilde{h}}(0,0)\right)\tilde{\delta h}$ on $\tilde{L}$ are invariant under the deck transformations, we have

\[
\|\delta h\|_{{\cal B}_1} = \|\tilde{\delta h}\|_{\tilde{{\cal B}}_1} \leq C \left\|\left(\frac{\partial \tilde{F}}{\partial \tilde{h}}(0,0)\right)\tilde{\delta h}\right\|_{\tilde{{\cal B}}_2} = C \left\|\left(\frac{\partial F}{\partial h}(0,0)\right)\delta h\right\|_{{\cal B}_2}.
\]
\hfill\rule{2.1mm}{2.1mm}\\
\begin{prop}
\label{fl}
\[
\left\|\frac{\partial F}{\partial h}(h,s)-\frac{\partial F}{\partial h}(0,0)\right\| \leq C\left(\nu_1^{\frac{1}{3}-\alpha} + [h]_{C^{2,\alpha}}\right).
\]
\[
\|F(0,s)\|_{{\cal B}_2} \leq C\nu_1^{\frac{1}{3}-\alpha}.
\]
\end{prop}
{\bf Proof:} Recall that
\[
\left(\frac{\partial F}{\partial h}(h,s)-\frac{\partial F}{\partial h}(0,0)\right) \delta h = (a_s^{ij} - a^{ij})\frac{\partial^2 \delta h}{\partial x_i\partial x_j} + (b_s^i - b^i)\frac{\partial\delta h}{\partial x_i}
\]

and
\[
F(0,s) = {\rm Im}(\log \eta_s(x,0)).
\]

Lemma \ref{fi} will imply the first estimate. The second estimate is implied by lemma \ref{fm}.
\hfill\rule{2.1mm}{2.1mm}\\
\begin{theorem}
\label{fd}
There exist $C_1,C_2>0$ such that for $L_{t,0}$ in a normal region, there exists a $C^1$-family of function $h_s\in U_{{\cal B}_1}(C_1\nu_1^{\frac{1}{3}-\alpha})$ such that $L_{t,s} = \phi_s(\{y = dh_s(x)\})$ is a smooth generalized special Lagrangian torus in $(X_{t,s},\tilde{\omega}_{t,s},\Omega_{t,s})$. Such $h_s$ is unique in $U_{{\cal B}_1}(C_2)$.
\end{theorem}
{\bf Proof:} According to proposition \ref{fb}, there exist $C_3>0$ such that

\[
\left\|\left(\frac{\partial F}{\partial h}(0,0)\right)^{-1}\right\| \leq C_3.
\]

By proposition \ref{fl}, there exist $C_1,C_2>0$ such that when $h \in U_{{\cal B}_1}(C_2)$ and $\nu_1$ is small, we have

\[
\left\|\frac{\partial F}{\partial h}(h,s)-\frac{\partial F}{\partial h}(0,0)\right\| \leq \frac{1}{2C_3}.
\]
\[
\|F(0,s)\|_{{\cal B}_2} \leq \frac{C_1\nu_1^{\frac{1}{3}-\alpha}}{2C_3}.
\]

Applying the implicit function theorem (theorem 3.2 in \cite{sl1}) to $F$, we have the conclusions of the theorem.
\hfill\rule{2.1mm}{2.1mm}\\

\se{Local generalized special Lagrangian fibration}
Although $L_{t,0}(r,c) = S_{r,c}\cap Y_t$ is Lagrangian in $(Y_t,\check{\omega}_t)$, since $\check{\omega}_t$ depends on $c$, when $(r,c)$ vary, $L_{t,0}(r,c)$ is not a Lagrangian fibration of $Y_t$ with respect to a fixed symplectic form independent of $(r,c)$. Notice that the symplectic deformation $\phi_s$ that deforms $L_{t,0}(r,c)$ to Lagrangian submanifold $\phi_1(L_{t,0}(r,c))$ in $(X_t=X_{t,1},\omega_t=\tilde{\omega}_{t,1})$ is also depending on $c$. Hence there is no obvious reason why when $(r,c)$ vary, $L_{t,1}(r,c)$ constructed in theorem \ref{fd} form a Lagrangian fibration of $(X_t=X_{t,1},\omega_t=\tilde{\omega}_{t,1})$. This problem will be clarified in this section. The argument can be separated into two parts. The first part is to argue that $(r,c) \rightarrow [L_{t,1}(r,c)]$ is a local embedding into the smooth local deformation space of generalized special Lagrangian torus in $(X_t=X_{t,1},\omega_t=\tilde{\omega}_{t,1})$. The second part is to argue that $L_{t,1}(r,c)$ do not intersect each other for different $(r,c)$. We will start with the first part.\\

Recall from section 4 that

\[
y_k = \lambda_k|z_k|^2 - \lambda_0|z_0|^2 - c_k,\ \ {\rm for}\ 1\leq k \leq l;
\]
\[
y_j = \rho_{j} - (\lambda_0|z_0|^2 - \eta) + \sum_{k\in I''}\frac{\lambda_{k,j}}{\lambda_k}(\lambda_k|z_k|^2 - \eta) - C^0_j,\ \ {\rm for}\ l+1\leq j \leq n.
\]

$x_k = \theta_k$ is a Darboux coordinate of $(Y_t, \hat{\omega}_t)$ such that $y=0$ when restricts to $L_{t,0}$. Lagrangian torus near $L_{t,0}$ can be parameterized by $H^1(L_{t,0},\mathbb{R})$. For a Lagrangian torus $L' = \{y(x)\}$ near $L_{t,0}$, the coordinate of such parameterization can be expressed as $[L']=\{\frac{1}{2\pi}\int y_k dx_k\}_{k=1}^n$. If $y$ is constant on $L'$, then $[L']=y$. Recall that

\[
\begin{CD}
(Y_t, \hat{\omega}_t) @>\varphi_c>> (Y_t,\check{\omega}_t) @>\phi_c>> (X_t,\omega_t)
\end{CD}
\]

are symplectomorphisms. (To stress the dependence on $c$, in this section, we are using $\varphi_c$, $\phi_c$ and $\psi_c = \phi_c\circ\varphi_c$ to replace $\varphi_1$, $\phi_1$ and $\psi_1 = \phi_1\circ\varphi_1$ in the previous notation.) Since $\omega_t$ is independent of $c$, $L'(r',c') = \psi_c^{-1}\circ\psi_{c'}(L_{t,0}(r',c'))$ is a Lagrangian torus in $(Y_t, \hat{\omega}_t(c))$ for $(r',c')$ near $(r,c)$. Notice that $L_{t,0}(r,c)$ is Hamiltonian equivalent to $\psi_c^{-1}(L_{t,1}(r,c))$. Through the symplectomorphism $\psi_c^{-1}$ the local deformation space of generalized special Lagrangian torus in $(X_t,\omega_t)$ can be identified with the local moduli space of Lagrangian torus (modulo Hamiltonian equivalence) near $L_{t,0}$ parameterized by $H^1(L_{t,0},\mathbb{R})$. Under such identification, $(r',c') \rightarrow [L_{t,1}(r',c')]$ can be reduced to $(r',c') \rightarrow [L'(r',c')]$. We will call $\displaystyle\left(\left\{\frac{c_k}{\nu_k^2}\right\}_{k=1}^l,\{\log r_k\}_{k=l+1}^n\right)$ the normalized coordinate of $(r,c)$. We will also normalize the coordinate on $H^1(L_{t,0},\mathbb{R})$, so that $[L']=\{\frac{1}{2\pi\nu_k^2}\int y_k dx_k\}_{k=1}^n$. A region in $\mathbb{R}^n$ will be called normal if for $(r,c)$ in it, $L_{t,0}(r,c)$ is in a normal region.\\
\begin{theorem}
\label{ga}
$(r,c) \rightarrow [L_{t,1}(r,c)]$ is a local embedding from a normal region in $\mathbb{R}^n$ into the smooth local deformation space of generalized special Lagrangian torus in $(X_t,\omega_t,\Omega_t)$. Under the normalized coordinates the tangent map and its inverse are both bounded in the normal region.
\end{theorem}
{\bf Proof:} We need to show that $(r',c') \rightarrow [L'(r',c')]$ is an embedding. Symbolically

\[
L'(r',c') - L_{t,0}(r,c) = (L_{t,0}(r',c') - L_{t,0}(r,c)) + (L'(r',c') - L_{t,0}(r',c')) 
\]
where
\[
L'(r',c') - L_{t,0}(r',c') = (\psi_c^{-1} - \psi_{c'}^{-1})\circ\psi_{c'}(L_{t,0}(r',c')).
\]

According to this, we can decompose $\delta y = \delta_1 y + \delta_2 y$. 

\[
\delta_1 y_k = \delta c_k,\ \ {\rm for}\ \ k\in I''.
\]

For $k\in I'$,

\[
\delta_1 y_k = \rho_{jk}\delta\log r_j^2 + \sum_{i\in I''}c_i(\log \lambda_i)_{jk}\delta\log r_j^2 + \sum_{j\in I''}(\log \lambda_j)_k\delta c_j + ((\log \Lambda)_k - 1)\delta \eta.
\]

Notice that

\[
\delta \eta = ((\log \Lambda)_j - 1 + {\rm Re}((\log p)_j)) \zeta \eta\delta\log r_j^2-\zeta\eta\sum_{j\in I''} \frac{\delta c_j}{c_j + \eta}.
\]

Consequently

\[
\delta_1 y_k = (\rho_{jk} + O(|z''|^2\rho_{jk})+ O(\eta))\delta\log r_j^2 + \sum_{j\in I''}O\left(\nu_k^2\nu_j^2 + \nu_0^2\right)\frac{\delta c_j}{c_j + \eta},\ \ {\rm for}\ \ k\in I'.
\]

\[
\delta_2 y_i = \frac{\partial y_i}{\partial \log z_j}\frac{\partial \log z_j}{\partial c_k}\delta c_k = \frac{\partial y_i}{\partial \log z_j}(c_k + \eta)\frac{\partial \log z_j}{\partial c_k}\frac{\delta c_k}{c_k + \eta}.
\]

Lemmas \ref{dn} and \ref{eg} imply that

\[
(c_k + \eta)\frac{\partial \log z_j}{\partial c_k} = O(A + B) = O(\epsilon^2).
\]

\[
\frac{\delta y_k}{\nu_k^2} = \frac{\delta c_k}{\nu_k^2} + O\left(\epsilon^2\left\{\frac{\delta c_i}{\nu_i^2}\right\}_{i\in I''}\right) ,\ \ {\rm for}\ \ k\in I''.
\]

\[
\frac{\delta y_k}{\nu_k^2} = \frac{\rho_{jk}}{\nu_k^2}\delta\log r_j^2 + O\left(\epsilon^2\left\{\delta\log r_i^2\right\}_{i\in I'}\right) + O\left(\epsilon^2\left\{\frac{\delta c_i}{\nu_i^2}\right\}_{i\in I''}\right),\ \ {\rm for}\ \ k\in I'.
\]

This expression of the tangent map under the normalized coordinates gives us the desired estimates.
\hfill\rule{2.1mm}{2.1mm}\\

For the second part of the argument, recall that the deformation of $L_{t,1}(r,c)$ in $(X_t=X_{t,1},\omega_t=\tilde{\omega}_{t,1})$ can be characterized by certain closed 1-forms on $L_{t,1}(r,c)$. To show that $L_{t,1}(r,c)$ do not intersect each other for different $(r,c)$, it is sufficient to show those closed 1-forms of deformation do not vanish anywhere on $L_{t,1}(r,c)$. It is straightforward to check that the deformation 1-forms of $L_{t,1}(r,c)$ are spanned by $\{dx_k - df_k\}_{k=1}^n$, where $\{f_k\}_{k=1}^n$ are functions on $L_{t,1}(r,c)$ satisfying

\[
a_1^{ij}\frac{\partial^2 f_k}{\partial x_i\partial x_j} + b_1^i\frac{\partial f_k}{\partial x_i} = b_1^k,
\]

which can be rewritten as

\[
a^{ij}\frac{\partial^2 f_k}{\partial x_i\partial x_j} + b^i\frac{\partial f_k}{\partial x_i} = b^k + (a^{ij} - a^{ij}_1)\frac{\partial^2 f_k}{\partial x_i\partial x_j} + (b^i - b^i_1)\frac{\partial f_k}{\partial x_i} + (b^k_1 - b^k).
\]

Lemmas \ref{fa}, \ref{fi} and propositions \ref{fb}, \ref{fl} imply that

\[
|f_k|_{C^{2,\alpha}} \leq C |b^k|_{C^{\alpha}}  + C\left(\nu_1^{\frac{1}{3}-\alpha} + [h]_{C^{2,\alpha}}\right)(|f_k|_{C^{2,\alpha}} + \frac{1}{\nu_k^2})
\]
\[
|f_k|_{C^{2,\alpha}} \leq C |b^k|_{C^{\alpha}} + \frac{C}{\nu_k^2}\left(\nu_1^{\frac{1}{3}-\alpha} + [h]_{C^{2,\alpha}}\right)\leq \frac{C\nu_1^{\frac{1}{3}-\alpha}}{\nu_k^2}.
\]

Applying lemma \ref{fg}, we have

\[
\left[\frac{\partial f_k}{\partial x_k} \right]_{C^0} \leq C\nu_k^{2 + \alpha}|f_k|_{C^{2,\alpha}} \leq C\nu_1^{\frac{1}{3}-\alpha}\nu_k^{\alpha}.
\]

Hence the $dx_k$-coefficient of $dx_k - df_k$ is of order $1+O((\nu_1)^{\frac{1}{3}-\alpha})$ which is non-zero when $\nu_1$ is small. Consequently $dx_k - df_k$ does not vanish anywhere.\\

Summing up the above argument and theorem \ref{ga}, we have\\
\begin{theorem}
\label{gb}
For $(r,c)$ in a normal region of $\mathbb{R}^n$, $L_{t,1}(r,c)$ constructed in theorem \ref{fd} form a generalized special Lagrangian torus fibration of an open set in $(X_t,\omega_t,\Omega_t)$.
\end{theorem}
\hfill\rule{2.1mm}{2.1mm}

\se{The uniqueness}
The uniqueness of generalized special Lagrangian submanifold under Hamiltonian deformation can be discussed from three perspectives --- the infinitesimal uniqueness, the local uniqueness and the global uniqueness. The infinitesimal uniqueness states that generalized special Lagrangian submanifold is isolated under Hamiltonian deformation, which is essentially implied by the work of Mclean (\cite{ML}). The global uniqueness attracted a lot of attention recently. In the cases that can be solved, one need among other things Floer (co)homology. For our application, we only need the local uniqueness, which states that when two Hamiltonian equivalent generalized special Lagrangian submanifolds are close enough to each other, they will coinside. Such local uniqueness can be proved essentially by the uniqueness part of our quantitative implicit function theorem (theorem 3.2 in \cite{sl1}). We will discuss this argument in this section.\\

{\bf Remark on notation:} The difference operator $\delta x = x_2 - x_1$ do not exactly behave as differential or derivative. Instead $\delta f(x) = f(x_2) - f(x_1) = \int_0^1 \frac{\partial f}{\partial x}(sx_2 + (1-s)x_1)\delta x ds$. The estimate of $\delta f(x)$ will come down to the estimate of $\frac{\partial f}{\partial x}(sx_2 + (1-s)x_1)\delta x$ for all $s\in [0,1]$. For this reason, we will use the slightly abused notation $\delta f(x) ``=" \frac{\partial f}{\partial x}\delta x$ to avoid unnecessary complication on notation.\\

Consider partitions $\{0,1,\cdots,n\} = (J'',J''',J') = (I'',I') = (\hat{I}'',\hat{I}')$, where $I''=J''$ and $\hat{I}' = J'$. $z=(z'',z''',z')$, $\tilde{z}=(\tilde{z}'',z''',z')$, $\tilde{z}'' = (z_0,z'')$. 

\[
S_{r,c} = \{\tilde{z} \in \mathbb{C}^{n+1}||z_j|=r_j, j\in I'; \lambda_k|z_k|^2 = c_k + \eta(z'), k\in I''.\}
\]

\[
\hat{S}_{\hat{r},\hat{c}} = \{\tilde{z} \in \mathbb{C}^{n+1}||z_j|=\hat{r}_j, j\in \hat{I}'; \hat{\lambda}_k|z_k|^2 = \hat{c}_k + \hat{\eta}(z'), k\in \hat{I}''.\}
\]

\[
X_{t,s} = \{z_0\cdots z_n = tp(s\tilde{z}'',z''',z')\},\ \ \hat{X}_{t,s} = \{z_0\cdots z_n = tp(s\tilde{z}'',sz''',z')\}.
\]

$Y_t = X_{t,0}$ and $\hat{Y}_t = \hat{X}_{t,0}$ are our local models and $X_{t,1} = \hat{X}_{t,1} = X_t$ is the actual hypersurface.\\ 

Let $L_{t,0} = S_{r,c}\cap X_{t,0}$, $\hat{L}_{t,0} = \hat{S}_{\hat{r},\hat{c}}\cap \hat{X}_{t,0}$. There are 4 coordinates $z^1$, $z^2 = \hat{\psi}_1^{-1}\circ\psi_1(z^1)$, $(x,y)^1$, $(x,y)^2$, involved in the argument. They can be viewed as coordinates on both $X_{t,0}$ and $\hat{X}_{t,0} = \hat{\psi}_1^{-1}\circ\psi_1(X_{t,0})$. To apply the implicit function theorem, we need to know the relation between $(L_{t,0},(x,y)^1)$ and $(\hat{L}_{t,1},(x,y)^1)$. We deduce this through the relation between $(L_{t,0},z^1)$ and $(\hat{L}_{t,1},z^1)$, which is achieved via the route

\[
(L_{t,0},z^1) \leftrightarrow (\hat{L}_{t,0},z^2) \leftrightarrow (\hat{L}_{t,0},z^1) \leftrightarrow (\hat{L}_{t,1},z^1).
\]

The relation between $(\hat{L}_{t,0},z^1) $ and $(\hat{L}_{t,1},z^1)$ is deduced through the relation between $(\hat{L}_{t,0},(x,y)^2)$ and $(\hat{L}_{t,1},(x,y)^2)$. (A word on terminology: When comparing $(L_{t,0},z^1) \leftrightarrow (\hat{L}_{t,0},z^2)$, we will identify $z^1$ and $z^2$, then compare $L_{t,0}$ and $\hat{L}_{t,0}$ under this common coordinate.)\\

The lemma \ref{df} can be rephrased as\\
\begin{lm}
\label{ha}
$\displaystyle\left(\frac{1}{\nu_j\nu_k}\frac{\partial y_j}{\partial \log |z_k|^2}\right)$ and $\displaystyle\left(\nu_j\nu_k\frac{\partial \log |z_j|^2}{\partial y_k}\right)$ are strongly T-bounded.\\

$\displaystyle\frac{\nu_j}{\nu_k}\frac{\partial x_j}{\partial \log |z_k|^2}$, $\displaystyle\frac{\nu_j}{\nu_k}\frac{\partial x_j}{\partial \theta_k}-\delta_{jk}$,$\displaystyle\frac{1}{\nu_j\nu_k}\frac{\partial y_j}{\partial \theta_k}$, $\displaystyle\nu_j\nu_k\frac{\partial \theta_j}{\partial y_k}$, $\displaystyle\frac{\nu_j}{\nu_k}\frac{\partial \theta_j}{\partial x_k}-\delta_{jk}$, $\displaystyle\frac{\nu_j}{\nu_k}\frac{\partial \log |z_j|^2}{\partial x_k}$ are of the order $\displaystyle O\left(\max\left(\frac{|z_0|^2}{|z_j|},\frac{|z_0|^2}{|z_k|}\right)\right)$.\\

Further more, the multi-derivatives of each non-constant term with respect to $\{\log z_k\}_{k=1}^n$ and their complex conjugates will hold same bound.
\end{lm}
\hfill\rule{2.1mm}{2.1mm}\\

The deformation $X_{t,0}\rightarrow X_{t,1} = X_t = \hat{X}_{t,1} \rightarrow \hat{X}_{t,0}$ results in coordinate transformation $z^1 \rightarrow z^2 = \hat{\psi}_1^{-1}\circ\psi_1(z^1)$. Applying lemmas \ref{dc} and \ref{ee} twice, we have
\begin{lm}
\label{he}
$\log z^1 - \log z^2$, $\displaystyle\frac{\partial \log z^1_j}{\partial \log z^2_k}-\delta_{jk}$, $\displaystyle\frac{\nu_j}{\nu_k}\frac{\partial \log z^1_j}{\partial \log z^2_k}-\delta_{jk}$  are of the order $\displaystyle O\left(\nu_1^{\frac{1}{3}}\right)$.\\

Further more, the multi-derivatives of each non-constant term with respect to $\{\log z_k\}_{k=1}^n$ and their complex conjugates will hold same bound.
\end{lm}
\hfill\rule{2.1mm}{2.1mm}\\

Being in the overlap of the two regions, we have $C\nu_0^{\frac{2}{3}} \leq \nu_j \leq C\nu_1^{\frac{2}{3}}$ for $j\in J'''$. Consequently, $C\nu_0^{\frac{2}{3}} \leq |z'''| \leq C\nu_1^{\frac{2}{3}}$.\\
\[
\hat{\lambda}_k = \lambda_k + O(|z'''|^2),\ \ p(s\tilde{z}'',z''',z') = p(s\tilde{z}'',sz''',z') + O(|z'''|).
\]

We have\\
\begin{lm}
\label{hb}
Assume that $\delta c_k = \hat{c}_k - c_k = O(\nu_k^2\nu_1^{\frac{1}{3}})$ for $k\in J''$, $\hat{c}_k - \hat{\lambda}_k r_k^2 = O(\nu_k^2\nu_1^{\frac{1}{3}})$ for $k\in J'''$ and $\delta \log r_k = \log \hat{r}_k - \log r_k = O(\nu_1^{\frac{1}{3}})$ for $k\in J'$, then between $(L_{t,0},z^1)$ and $(\hat{L}_{t,0},z^2)$, we have
\[
\delta \log \eta = O(\nu_1^{\frac{1}{3}}),\ \ \delta \log |z_k| = O(\nu_1^{\frac{1}{3}}),
\]
and for $\beta = 0,\alpha$
\[
[\delta \log \eta]_{C^{1,\beta}} = O(\nu_1^{-\beta}),\ \ [\delta \log |z_k|]_{C^{1,\beta}} = O(\nu_k^{-2}\nu_1^{2-\beta}).
\]
\end{lm}
{\bf Proof:} According to the assumptions of the lemma

\[
\delta \log |z_k| = \delta \log r_k = O(\nu_1^{\frac{1}{3}}),\ \ {\rm for}\ k\in J'.
\]
\[
\delta \log |z_k|^2 = \log \left(\frac{\hat{c}_k + \eta}{\hat{\lambda}_k r_k^2}\right) = O\left(\frac{\hat{c}_k - \hat{\lambda}_k r_k^2 + \eta}{\hat{\lambda}_k r_k^2}\right) = O(\nu_1^{\frac{1}{3}}),\ \ {\rm for}\ k\in J'''.
\]

For $k\in J''$

\[
\delta \log |z_k|^2 = \frac{\delta \eta + \delta c_k}{\eta + c_k} - \delta \log \lambda_k.
\]
\[
\delta \log \lambda_k = \delta \log \lambda_k(0,0,z') + (\log \lambda_k(0,z''',z') - \log \lambda_k(0,0,z')) 
\]
\[
= \sum_{j\in J'}2\frac{\lambda_{k,j}}{\lambda_k}(0,0,z')\frac{\delta r_j}{r_j} + O(|z'''|^2) = O(\nu_1^{\frac{1}{3}}).
\]
\[
\delta \log |p|^2 = \delta \log |p|^2(0,0,z') + (\log |p|^2(0,z''',z') - \log |p|^2(0,0,z')) 
\]
\[
= \sum_{k\in J'}2{\rm Re}\left(\frac{p_k}{p}\right)(0,0,z')\frac{\delta r_k}{r_k} + O(|z'''|) = O(\nu_1^{\frac{1}{3}}).
\]

Substitute all these to $\displaystyle\sum_{k=0}^n \delta \log |z_k|^2 = \delta \log |p|^2$, we have

\[
\zeta_1^{-1}\delta \log \eta + \sum_{k\in J''}\left(\frac{\delta c_k}{\eta + c_k}\right) = O(\nu_1^{\frac{1}{3}}).
\]

Consequently $\delta \log \eta = O(\nu_1^{\frac{1}{3}})$ and $\delta \log |z_k|^2 = O(\nu_1^{\frac{1}{3}})$, for all $k$.\\

Notice that only $\eta$ and $p$ are depending on $\theta$. We have

\[
\frac{\partial \log \eta}{\partial \theta_k} = \zeta \frac{\partial \log |p|^2}{\partial \theta_k} = O(\nu_k).
\]
\[
\frac{\partial \delta \log |z_j|}{\partial \theta} = 0,\ \ {\rm for}\ j\in J'.
\]
\[
\frac{\partial \delta \log |z_j|^2}{\partial \theta_k} = \frac{\eta}{\hat{c}_j + \eta}\frac{\partial \log \eta}{\partial \theta_k} = O(\nu_k\nu_j^{-2}\nu_1^2),\ \ {\rm for}\ j\in J'''.
\]
\[
\frac{\partial \delta \log |z_j|^2}{\partial \theta_k} = \frac{\eta}{c_j + \eta}\frac{\partial \delta \log \eta}{\partial \theta_k} + \frac{c_j\eta\delta \log \eta}{(c_j + \eta)^2}\frac{\partial \log \eta}{\partial \theta_k},\ \ {\rm for}\ j\in J''.
\]
\[
\frac{\partial \delta \log |p|^2}{\partial \theta_k} = O(\nu_k\nu_1^{\frac{1}{3}})\ \ {\rm for}\ k\not\in J'''.
\]
\[
\frac{\partial \delta \log |p|^2}{\partial \theta_k} = O(\nu_k)\ \ {\rm for}\ k\in J'''.
\]

Substitute all these to $\displaystyle\sum_{j=0}^n \frac{\partial \delta \log |z_j|^2}{\partial \theta_k} = \frac{\partial \delta \log |p|^2}{\partial \theta_k}$, we get $\displaystyle\zeta_1\frac{\partial \delta \log \eta}{\partial \theta_k} = O(\nu_k)$. Consequently

\[
\frac{\partial \delta \log \eta}{\partial \theta_k} = O(\nu_k),\ \ \frac{\partial \delta \log |z_j|^2}{\partial \theta_k} = O(\nu_k\nu_j^{-2}\nu_1^2)
\]

Further more, it is straightforward to derive that additional derivatives with respect to $\theta$ will hold the same bound. Consequently

\[
[\delta \log \eta]_{C^1} = O(1),\ \ [\delta \log |z_k|]_{C^1} = O(\nu_k^{-2}\nu_1^{2}).
\]

Notice that $\nu_1$ is the smallest scale of the torus. For any function $f(\theta)$, it is easy to verify that

\[
[f]_{C^\alpha} \leq C\nu_1^{-\alpha}\left[\frac{\partial f}{\partial \theta}\right]_{C^0}\leq C\nu_1^{-\alpha}\left[\frac{\partial f}{\partial \theta}\right]_{C^0}.
\]

Hence

\[
[\delta \log \eta]_{C^{1,\alpha}} = O(\nu_1^{-\alpha}),\ \ [\delta \log |z_k|]_{C^{1,\alpha}} = O(\nu_k^{-2}\nu_1^{2-\alpha}).
\]
\hfill\rule{2.1mm}{2.1mm}\\

{\bf Remark:} In every lemma of this section, the $C^{1,\alpha}$-estimate will always be related to the $C^1$-estimate in the same way as in lemma \ref{hb}. Therefore, from now on, we will only discuss the $C^1$-estimates and omit the discussion of the $C^{1,\alpha}$-estimates.\\
\begin{lm}
\label{hf}
Between $(\hat{L}_{t,0},z^1)$ and $(\hat{L}_{t,0},z^2)$, we have
\[
\delta \log |z_k| = O(\nu_1^{\frac{1}{3}}),\ \ {\rm and}\ \ [\nu_k\delta \log |z_k|]_{C^{1,\beta}} = O(\nu_1^{\frac{1}{3}-\beta}),
\]
for $\beta = 0,\alpha$.
\end{lm}
{\bf Proof:} The equation of $\hat{L}_{t,0}$ is $(\log |z^2|(\theta^2),\theta^2)$ under $z^2$. The equation of $\hat{L}_{t,0}$ is $(\log |z^1|(\log |z^2|(\theta^2),\theta^2),\theta^1(\log |z^2|(\theta^2),\theta^2))$ under $z^1$. Let $\theta^2 = \Psi(\theta^1)$ be the inverse function of $\theta^1 = \theta^1(\log |z^2|(\theta^2),\theta^2)$. By lemma \ref{he}

\begin{equation}
\label{hm}
\frac{\nu_j}{\nu_k}\frac{\partial \Psi_j}{\partial \theta^1_k} = \frac{\nu_j}{\nu_k}\left(\frac{\partial \theta^1}{\partial \theta^2} + \frac{\partial \theta^1}{\partial \log |z^2|}\frac{\partial \log |z^2|}{\partial \theta^2}\right)^{-1}_{jk} = I + O(\nu_1^{\frac{1}{3}}).
\end{equation}

The equation of $\hat{L}_{t,0}$ is $(\log |z^1|(\log |z^2|(\Psi(\theta^1)),\Psi(\theta^1)),\theta^1)$ under $z^1$ after the change of parameter $\theta^2 = \Psi(\theta^1)$. When we compare $(\hat{L}_{t,0},z^1)$ and $(\hat{L}_{t,0},z^2)$, we identify $z_1$ and $z_2$. In particular, let $\theta = \theta_1 = \theta_2$ and $\tilde{\theta} = \Psi(\theta)$. Then

\[
\delta \log |z|(\theta) = \log |z^2|(\theta) - \log |z^1|(\log |z^2|(\tilde{\theta}),\tilde{\theta}),
\]

\[
\delta \log |z|(\theta) = \log |z^2|(\theta) - \log |z^2|(\tilde{\theta}) + \log |z^2|(\tilde{\theta}) - \log |z^1|(\log |z^2|(\tilde{\theta}), \tilde{\theta}),
\]
\[
= (\theta - \tilde{\theta})\frac{\partial \log |z^2|}{\partial  \theta^2} + (\log |z^2| - \log |z^1|).
\]

\[
\frac{\partial \delta \log |z|}{\partial \theta} = \frac{\partial \log |z^2|}{\partial \theta^2}(\theta) - \left(\frac{\partial \log |z^1|}{\partial \theta^2}(\tilde{\theta}) + \frac{\partial \log |z^1|}{\partial \log |z^2|}\frac{\partial \log |z^2|}{\partial \theta^2}(\tilde{\theta})\right)\frac{\partial \Psi}{\partial \theta}
\]

Applying lemma \ref{he} and (\ref{hm}), we have $\delta \log |z|(\theta) = O(\nu_1^{\frac{1}{3}})$ and 

\[
\frac{\nu_j}{\nu_k}\frac{\partial \delta \log |z_j|}{\partial \theta_k} = \frac{\nu_j}{\nu_k}\left(\frac{\partial \log |z^2_j|}{\partial \theta^2_k}(\theta) - \frac{\partial \log |z^2_j|}{\partial \theta^2_k}(\tilde{\theta})\right) + O(\nu_1^{\frac{1}{3}}) = O(\nu_1^{\frac{1}{3}})
\] 
\hfill\rule{2.1mm}{2.1mm}
\begin{lm}
\label{hd}
Assume that between $(L_{t,0},z^1)$ and $(\hat{L},z^1)$, we have
\[
\delta \log |z_k| = O(\nu_1^{\frac{1}{3}-\alpha_1}),\ \ {\rm and}\ \ [\nu_k\delta \log |z_k|]_{C^{1,\beta}} = O(\nu_1^{\frac{1}{3}-\alpha_1-\beta}),
\]
for $\beta = 0,\alpha$. Then between $(L_{t,0},(x,y)^1)$ and $(\hat{L},(x,y)^1)$, we have

\[
y_k = O(\nu_k^2\nu_1^{\frac{1}{3}-\alpha_1}),\ \ {\rm and}\ \ [y_k]_{C^{1,\beta}} = O(\nu_k\nu_1^{\frac{1}{3}-\alpha_1-\beta}),
\]
for $\beta = 0,\alpha$.
\end{lm}
{\bf Proof:} Recall that $y|_{L_{t,0}}=0$. By assumption and lemma \ref{ha}, we have
\[
\delta y_k = y_k = \frac{\partial y_k}{\partial \log |z|}\delta \log |z| = O(\nu_k^2\nu_1^{\frac{1}{3}-\alpha_1}).
\]
\[
\frac{\partial y}{\partial x} = \left(\left(\frac{\partial}{\partial \theta}\left(\frac{\partial y}{\partial \log |z|}\right) + \frac{\partial}{\partial \log |z|}\left(\frac{\partial y}{\partial \log |z|}\right)\frac{\partial \log |z|}{\partial \theta}\right)\delta \log |z| \right.
\]
\[
\left.+ \frac{\partial y}{\partial \log |z|}\frac{\partial \delta \log |z|}{\partial \theta}\right)\left(\frac{\partial x}{\partial \theta} + \frac{\partial x}{\partial \log |z|}\frac{\partial \log |z|}{\partial \theta}\right)^{-1}.
\]

Lemma \ref{ha} implies that 

\[
\left[\frac{1}{\nu_j\nu_k}\frac{\partial y_j}{\partial x_k}\right]_{C^0} \leq C([\delta \log |z|]_{C^0} + [\nu\delta \log |z|]_{C^1}) = O(\nu_1^{\frac{1}{3}-\alpha_1}).
\]
\hfill\rule{2.1mm}{2.1mm}
\begin{lm}
\label{hc}
Under the same assumption as in lemma \ref{hb}, between $(L_{t,0},z^1)$ and $(\hat{L}_{t,1},z^1)$, we have
\[
\delta \log |z_k| = O(\nu_1^{\frac{1}{3}-\alpha}),\ \ {\rm and}\ \ [\nu_k\delta \log |z_k|]_{C^{1,\beta}} = O(\nu_1^{\frac{1}{3}-\alpha-\beta}),
\]
for $\beta = 0,\alpha$.
\end{lm}
{\bf Proof:} According to lemmas \ref{hb} and \ref{hf}, we only need to prove the estimate between $(\hat{L}_{t,0},z^1)$ and $(\hat{L}_{t,1},z^1)$, which can be reduced to the estimate between $(\hat{L}_{t,0},(x,y)^2)$ and $(\hat{L}_{t,1},(x,y)^2)$. Let $z(x,y)$ be the composition $(x^2,y^2)\rightarrow z^2 \rightarrow z^1$. Lemmas \ref{ha} and \ref{he} together give the estimates for $z(x,y)$.

\begin{equation}
\label{hg}
\delta \log |z_k|(\theta) = \log |z_k|(x(\theta),y(x(\theta))) - \log |z_k|(x'(\theta),0),
\end{equation}

where $x(\theta)$ and $x'(\theta)$ are the inverse functions of $\theta = \theta(x,y(x))$ and $\theta = \theta(x',0)$.

\[
\delta \log |z_k| = \frac{\partial \log |z_k|}{\partial x_j}\delta x_j + \frac{\partial \log |z_k|}{\partial y_j}y_j,
\]
\[
0 = \delta \theta_k = \frac{\partial \theta_k}{\partial x_j}\delta x_j + \frac{\partial \theta_k}{\partial y_j}y_j,
\]
\begin{equation}
\label{hn}
\delta \log |z_k| = \left(\frac{\partial \log |z_k|}{\partial y_j}-\frac{\partial \log |z_k|}{\partial x}\left(\frac{\partial \theta}{\partial x}\right)^{-1}\frac{\partial \theta}{\partial y_j}\right) y_j.
\end{equation}

Lemmas \ref{ha} and \ref{he} imply that 

\begin{equation}
\label{hh}
\nu_j^2\frac{\partial \log |z|}{\partial y_j},\ \ \frac{\partial \log |z|}{\partial x},\ \ \left(\frac{\partial \theta}{\partial x}\right)^{-1},\ \ \nu_j^2\frac{\partial \theta}{\partial y_j}
\end{equation}

are all bounded. Also recall that $y(x) = dh(x)$. Applying lemma \ref{fg} and theorem \ref{fd}, we have

\[
[\delta \log |z_k|]_{C^0} \leq C\max_j\left(\frac{y_j}{\nu_j^2}\right) \leq C[h]_{C^2} = O(\nu_1^{\frac{1}{3}-\alpha}).
\]

Using (\ref{hg}), it is straightforward to derive

\begin{equation}
\label{hi}
\frac{\partial \delta\log |z|}{\partial \theta}\! =\! \delta\!\left(\!\frac{\partial \log |z|}{\partial x}\!\left(\!\frac{\partial \theta}{\partial x}\!\right)^{-1}\right)\! +\! \left(\!\frac{\partial \log |z|}{\partial y}\! -\! \frac{\partial \log |z|}{\partial x}\left(\!\frac{\partial \theta}{\partial x}\!\right)^{-1}\!\frac{\partial \theta}{\partial y}\!\right)\!\frac{\partial y}{\partial \theta}.
\end{equation}

Here we used the formulas

\[
\frac{\partial x}{\partial \theta} = \left(\frac{\partial \theta}{\partial x} + \frac{\partial \theta}{\partial y}\frac{\partial y}{\partial x}\right)^{-1},\ \ \frac{\partial x'}{\partial \theta} = \left(\frac{\partial \theta}{\partial x'}\right)^{-1},
\]
\[
\left(\frac{\partial \theta}{\partial x} + \frac{\partial \theta}{\partial y}\frac{\partial y}{\partial x}\right)^{-1} - \left(\frac{\partial \theta}{\partial x}\right)^{-1} = - \left(\frac{\partial \theta}{\partial x}\right)^{-1}\frac{\partial \theta}{\partial y}\frac{\partial y}{\partial \theta}.
\]

Similar as in (\ref{hn}), we have

\[
\delta\left(\frac{\partial \log |z|}{\partial x}\left(\frac{\partial \theta}{\partial x}\right)^{-1}\right) = \left(\frac{\partial}{\partial y_j}\left(\frac{\partial \log |z|}{\partial x}\left(\frac{\partial \theta}{\partial x}\right)^{-1}\right)\right.
\]
\[
\left.-\frac{\partial}{\partial x}\left(\frac{\partial \log |z|}{\partial x}\left(\frac{\partial \theta}{\partial x}\right)^{-1}\right)\left(\frac{\partial \theta}{\partial x}\right)^{-1}\frac{\partial \theta}{\partial y_j}\right) y_j.
\]

Lemmas \ref{ha} and \ref{he} imply that $\displaystyle\frac{\nu_j}{\nu_k}\left(\frac{\partial \log |z|}{\partial x}\left(\frac{\partial \theta}{\partial x}\right)^{-1}\right)_{jk}$ and their multi-derivatives with respect to $\log z$ are bounded. Apply lemmas \ref{ha} and \ref{he}, we get that terms as in (\ref{hh}) replacing $\log |z|$ by $\displaystyle\frac{\nu_j}{\nu_k}\left(\frac{\partial \log |z|}{\partial x}\left(\frac{\partial \theta}{\partial x}\right)^{-1}\right)_{jk}$ are bounded, which implies

\begin{equation}
\label{hj}
\left[\frac{\nu_j}{\nu_k}\delta\left(\frac{\partial \log |z|}{\partial x}\left(\frac{\partial \theta}{\partial x}\right)^{-1}\right)_{jk}\right]_{C^0} \leq C[h]_{C^2}.
\end{equation}

Notice that

\[
\frac{\partial y}{\partial \theta} = \frac{\partial y}{\partial x}\frac{\partial x}{\partial \theta} = \frac{\partial y}{\partial x}\left(\frac{\partial \theta}{\partial x} - \frac{\partial \theta}{\partial y}\frac{\partial y}{\partial x}\right)^{-1}.
\]

Lemmas \ref{ha} and \ref{he} imply that 

\[
\nu_j\nu_k\frac{\partial \log |z_k|}{\partial y_j},\ \ \frac{\nu_j}{\nu_k}\frac{\partial \log |z_j|}{\partial x_k},\ \ \frac{\nu_j}{\nu_k}\left(\frac{\partial \theta}{\partial x}\right)^{-1}_{jk},\ \ \frac{1}{\nu_j\nu_k}\frac{\partial \theta_k}{\partial y_j}
\]

are all bounded. Consequently

\begin{equation}
\label{hk}
\left[\frac{\nu_j}{\nu_k}\left(\frac{\partial \log |z_j|}{\partial y} - \frac{\partial \log |z|}{\partial x}\left(\frac{\partial \theta}{\partial x}\right)^{-1}\frac{\partial \theta}{\partial y}\right)\frac{\partial y}{\partial \theta_k}\right]_{C^0} \leq C[h]_{C^2}.
\end{equation}

Combine (\ref{hj}) and (\ref{hk}) together with (\ref{hi}), we have

\[
[\nu_j\delta \log |z_j|]_{C^1} \leq C[h]_{C^2} = O(\nu_1^{\frac{1}{3}-\alpha}).
\]
\hfill\rule{2.1mm}{2.1mm}
\begin{theorem}
\label{hl}
For suitable $(r,c)$ and $(\hat{r},\hat{c})$, the generalized special Lagrangian torus $L_{t,1}(r,c)$ coincides with $\hat{L}_{t,1}(\hat{r},\hat{c})$.
\end{theorem}
{\bf Proof:} First choose $(r,c)$ and $(\hat{r},\hat{c})$ so that $\delta c_k = 0$ for $k\in J''$, $\hat{c}_k = \hat{\lambda}_k r_k^2$ for $k\in J'''$ and $\delta \log r_k = 0$ for $k\in J'$. Then lemmas \ref{hb} and \ref{hf} imply that between $(L_{t,0},z^1)$ and $(\hat{L}_{t,0},z^1)$

\[
\delta \log |z_k| = O(\nu_1^{\frac{1}{3}}),\ \ {\rm and}\ \ [\nu_k\delta \log |z_k|]_{C^{1,\beta}} = O(\nu_1^{\frac{1}{3}-\beta}),
\]
for $\beta = 0,\alpha$. Applying lemma \ref{hd} for $\alpha_1=0$ and $\hat{L}=\hat{L}_{t,0}$, we see that between $(L_{t,0},(x,y)^1)$ and $(\hat{L}_{t,0},(x,y)^1)$, we have

\[
y_k = O(\nu_k^2\nu_1^{\frac{1}{3}}),\ \ {\rm and}\ \ [y_k]_{C^{1,\beta}} = O(\nu_k\nu_1^{\frac{1}{3}-\beta}),
\]
for $\beta = 0,\alpha$. In particular, $\hat{L}_{t,0}$ as an element of the local moduli space of Lagrangian torus (modulo Hamiltonian equivalence) near $L_{t,0}$ parameterized by $H^1(L_{t,0},\mathbb{R})$ has coordinate $\{\frac{1}{2\pi\nu_k^2}\int y_k dx_k\}_{k=1}^n = O(\nu_1^{\frac{1}{3}})$ under the normalized coordinate of $H^1(L_{t,0},\mathbb{R})$. Theorem \ref{ga} implies that one may adjust $(\hat{r},\hat{c})$ so that $(\hat{L}_{t,0},(x,y)^1)$ is harmiltonian equivalent to $(L_{t,0},(x,y)^1)$ and the assumption of lemma \ref{hb} is still satisfied. Then lemma \ref{hc} implies that between $(L_{t,0},z^1)$ and $(\hat{L}_{t,1},z^1)$

\[
\delta \log |z_k| = O(\nu_1^{\frac{1}{3}-\alpha}),\ \ {\rm and}\ \ [\nu_k\delta \log |z_k|]_{C^{1,\beta}} = O(\nu_1^{\frac{1}{3}-\alpha-\beta}),
\]
for $\beta = 0,\alpha$. Applying lemma \ref{hd} for $\alpha_1=\alpha$ and $\hat{L}=\hat{L}_{t,1}$, we see that between $(L_{t,0},(x,y)^1)$ and $(\hat{L}_{t,1},(x,y)^1)$, we have

\[
y_k = O(\nu_k^2\nu_1^{\frac{1}{3}-\alpha}),\ \ {\rm and}\ \ [y_k]_{C^{1,\beta}} = O(\nu_k\nu_1^{\frac{1}{3}-\alpha-\beta}),
\]

for $\beta = 0,\alpha$. Since $\hat{L}_{t,1}$ is Hamiltonian equivalent to $\hat{L}_{t,0}$, $\hat{L}_{t,1}$ is also Hamiltonian equivalent to $L_{t,0}$. Namely, under $(x,y)^1$, $\hat{L}_{t,1}$ can be characterized as $y = d\hat{h}$. and $|\hat{h}|_{C^{2,\alpha}} = O(\nu_1^{\frac{1}{3}-2\alpha})$. When $\frac{1}{3}-2\alpha>0$ and $\nu_1$ is small, we have $|\hat{h}|_{C^{2,\alpha}}< C_2$ for $C_2$ in theorem \ref{fd}. Applying the uniqueness part of theorem \ref{fd}, $L_{t,1}(r,c)$ coincides with $\hat{L}_{t,1}(\hat{r},\hat{c})$.
\hfill\rule{2.1mm}{2.1mm}\\

The case of $I'' = \{0\}$ need some additional comments. In this case, since $1\not\in I''$, $\nu_1$ does not play the same role as in the other cases, where $|I''|\geq 2$. If one traces through all our argument, it is easy to observe that all the estimates and results before section 8 are still valid in this case as long as one replace $\nu_1$ by $\nu_0$. Since $\nu_0\leq \nu_1$, clearly all the results thus far are still valid when $I'' = \{0\}$.\\

{\bf Remark:} The case of $I'' = \{0\}$ corresponds to fibrations over top dimensional faces of $\partial \Delta$. If one works through the estimates and results before section 8 more carefully, one can actually show that the fibration can be constructed without assuming $|p(z)| \geq C >0$ in the case of $I'' = \{0\}$. Instead of exploring this argument further, in the next section, we will use an alternative argument by making use of the fibration over top dimensional faces of $\partial \Delta$ we constructed earlier in theorems 5.1 and 5.2 of \cite{sl1} and more quantitatively in theorem 4.1 of \cite{sl2} using a different method.\\

\se{Monodromy representing generalized special Lagrangian torus fibrations}
Let $F_\omega: P_\Delta \rightarrow \Delta \subset M_{\mathbb{R}}$ be the moment map determined by a toric metric $\omega$ on the toric variety $P_\Delta$. We assume that the fan $\Sigma$ of $(n+1)$-dimensional $P_\Delta$ consists only simplicial cones. Then $P_\Delta$ will have at most finite toric orbifold singularities. We also assume that $\omega$ is a smooth toric orbifold metric on $P_\Delta$. (Any toric variety will always has a crepant resolution being such simplicial toric variety.) The set of vertices of $\Delta$ can be identified to the set $\Sigma(n+1)$ of top dimensional cones in $\Sigma$. We will use $m_\sigma$ to denote the vertex of $\Delta$ corresponding to $\sigma \in \Sigma(n+1)$. Using the bari-center subdivision of $\Delta$, we have the following decomposition of $\Delta$ by convex polyhedrons.

\[
\Delta = \bigcup_{\sigma\in \Sigma(n+1)} \Delta_\sigma,
\]

where $\Delta_\sigma$ is a polyhedron neighborhood of the vertex $m_\sigma$ of $\Delta$. (More precisely, $\Delta_\sigma$ is the union of closed simplices in the bari-center subdivision of $\Delta$ containing $m_\sigma$.) Let $\tilde{\Delta}_\sigma$ be an open neighborhood of $m_\sigma$ that is slightly larger than $\Delta_\sigma$. Then $F_\omega^{-1}(\tilde{\Delta}_\sigma)$ is an open set in the affine toric piece ${\rm Spec}(\mathbb{C}[\sigma^\vee])$. Since $\sigma$ is simplicial, there exists a finite orbifold cover $\pi_\sigma: \mathbb{C}^{n+1} \rightarrow {\rm Spec}(\mathbb{C}[\sigma^\vee]) \cong \mathbb{C}^{n+1}/\Pi$, where $\Pi$ is a finite abelian group of toric isometries of $\mathbb{C}^{n+1}$. Let $\tilde{F}_\omega = F_\omega \circ \pi_\sigma: \mathbb{C}^{n+1} \rightarrow \Delta$. There exist constants $C_1,C_2>0$ such that

\[
\{|z|< C_1\} \subset \tilde{F}_\omega^{-1}(\tilde{\Delta}_\sigma) \subset \{|z|< C_2\}.
\]

For each decomposition $\{0,\cdots,n\} = I'' \cup I'$ with $|I''|\geq 1$, let $\tilde{\nu}_0$ and $\tilde{\nu}_1$ denote the smallest and the second smallest elements in $\{|z_i||i\in I''\}$. (If $|I''|=1$ we take $\tilde{\nu}_0=\tilde{\nu}_1$.) We may define the following regions

\[
\tilde{U}_{\sigma,I'} = \left\{z\in \mathbb{C}^{n+1}\left|\begin{array}{l}|z| <C_2,\ |p(z)|\geq C_3,\ |z_i| \leq |z_j|,\\|z_i| \leq C_4\tilde{\nu}_1^{\frac{2}{3}},\ |z_j| \geq C_4\tilde{\nu}_0^{\frac{2}{3}},\ {\rm for}\ i\in I'', j\in I'.\end{array}\right.\right\}
\]

and $U_{\sigma,I'} = \pi_\sigma (\tilde{U}_{\sigma,I'})$.\\

Recall that the family of Calabi-Yau hypersurfaces in $P_\Delta$ is defined as

\[
X_t = \{\tilde{s}_t^{-1}(0)\},\  {\rm where}\ \tilde{s}_t = s_{m_o} + ts,\ \ s = \sum_{m\in \Delta_0\setminus \{m_o\}} a_m s_m.
\]

$\pi_\sigma^{-1}(X_t) \subset \mathbb{C}^{n+1}$ is defined by $z_0 \cdots z_n = tp(z)$, where $p(z) = (s/s_{m_\sigma})\circ \pi_\sigma(z)$. $z_0 \cdots z_n$ and toric monomials in $p(z)$ are invariant under $\Pi$-action. Let

\[
\tilde{U}_{\sigma} = \left\{z\in \tilde{F}_\omega^{-1}(\tilde{\Delta}_\sigma):|p(z)|\geq C_3\right\}
\]

and $U_{\sigma} = \pi_\sigma (\tilde{U}_{\sigma})$. Then clearly

\[
\tilde{U}_{\sigma} \subset \bigcup_{I'} \tilde{U}_{\sigma,I'}.
\]

\begin{theorem}
\label{ia}
For $t$ small enough, there exists a generalized special Lagrangian torus fibration over $\partial \Delta$ outside a small neighborhood of $F_\omega(X_t \cap {\rm Sing}(X_0))$, which is the singular locus of the Lagrangian torus fibration constructed in \cite{toric} via the Hamiltonian-gradient flow. Such (partial) fibration is called monodromy representing in \cite{sl1,sl2}.
\end{theorem}
{\bf Proof:} It is easy to see that $\tilde{U}_{\sigma,I'} \subset \mathbb{C}^{n+1}$ near $\pi_\sigma^{-1}(X_t)$ is a normal region as defined in the beginning of section 6. Theorem \ref{gb} implies that we can construct generalized special Lagrangian torus fibration for $\tilde{U}_{\sigma,I'} \cap \pi_\sigma^{-1}(X_t)$. Theorem \ref{hl} implies that such fibrations for different decompositions $I''\cup I'$ will coincide in the overlap and piece together to form the generalized special Lagrangian torus fibration for $\tilde{U}_{\sigma} \cap \pi_\sigma^{-1}(X_t)$.\\

It is straightforward to check that all objects and processes of our construction are invariant under $\Pi$-action. Consequently, the resulting fibration is invariant under $\Pi$-action and decends to form the generalized special Lagrangian torus fibration for $U_{\sigma} \cap X_t$. Such fibrations for different $\sigma\in \Sigma(n+1)$ in the overlap are constructed with the same local model and deformation process. They will coincide and piece together to form the generalized special Lagrangian torus fibration for $X_t$ outside a small neighborhood of $F_\omega^{-1}(F_\omega(X_t \cap X_0))$.\\

Glue in the generalized special Lagrangian fibration over top dimensional faces of $\partial \Delta$ constructed in theorem 4.1 of \cite{sl2}, we get our monodromy representing generalized special Lagrangian torus fibration for $X_t$ outside a small neighborhood of $F_\omega^{-1}(F_\omega(X_t \cap {\rm Sing}(X_0)))$. (Recall that $X_t \cap {\rm Sing}(X_0)$ and $F_\omega(X_t \cap {\rm Sing}(X_0))$ are the singular set and the singular locus of the Lagrangian fibration constructed in \cite{toric} via the Hamiltonian-gradient flow.) The smallness of the neighborhood is determined by how small $C_3$ is, which can be made arbitrarily small by taking $t$ small.\\

The argument for the final gluing is essentially the same as the arguments used in theorem 7.1 of \cite{sl1} and theorem 5.1 of \cite{sl2}. One only need to show that the two fibrations coincide at a torus fibre of bounded geometry over the top dimensional face and away from $X_t \cap X_0$. The two local models are 

\[
L_{t,0} = \{\tilde{z}\in X_{t,0}||z_k|=r_k ({\rm constant}),\ {\rm for}\ 1\leq k \leq n\}
\]
\[
L_0 = \{z_0=0, |z_k|=r_k ({\rm constant}),\ {\rm for}\ 1\leq k \leq n\}
\]

where $\log r_k = O(1)$ for $1\leq k \leq n$. Both $L_{t,0}$ and $L_0$ are of bounded geometry. It is easy to check explicitly that $\|L_{t,0} - L_0\|_{C^{1,\alpha}} = O(t)$. $\|L_t - L_0\|_{C^{1,\alpha}} = O(t)$ for the same reason as in the proof of theorem 7.1 in \cite{sl1}. Since $(X_{t,s}, \tilde{\omega}_{t,s})$ are $O(t)$-perturbations of $(X_{t,0}, \tilde{\omega}_{t,0})$, we also have $\|L_{t,1} - L_{t,0}\|_{C^{1,\alpha}} = O(t)$ for the same reason as in the proof of theorem 7.1 in \cite{sl1}. Consequently $\|L_{t,1} - L_t\|_{C^{1,\alpha}} = O(t)$. Then the rest of the argument in the proof of theorem 7.1 in \cite{sl1} gives us the desired gluing.
\hfill\rule{2.1mm}{2.1mm}\\

\ifx\undefined\bysame
\newcommand{\bysame}{\leavevmode\hbox to3em{\hrulefill}\,}
\fi

\end{document}